\documentclass{gtart}

\def\ifplaintex{\expandafter\ifx\csname documentclass\endcsname\relax}

\def\gtp{{\mathsurround=0pt\it $\cal G\mskip-2mu$eometry \&\ 
$\cal T\!\!$opology $\cal P\!$ublications}}  

\def\recd{{\small Received:\qua\receiveddate\ifx\reviseddate\relax
\else\qquad Revised:\qua\reviseddate\fi\par}} 

\def\lognumber#1{\def\thelognumber{#1}}
\def\volumenumber#1{\def\thevolumenumber{#1}}
\def\volumeyear#1{\def\thevolumeyear{#1}}
\def\papernumber#1{\def\thepapernumber{#1}}
\def\pagenumbers#1#2{\def\startpage{#1}\def\finishpage{#2}}
\def\published#1{\def\publishdate{#1}}

\def\received#1{\def\receiveddate{#1}}
\def\revised#1{\def\reviseddate{#1}}
\def\accepted#1{\def\accepteddate{#1}}

\def\asciiaddress#1{\def\theasciiaddress{#1}}

\long\def\asciiabstract#1{\long\def\theasciiabstract{#1}}

\let\\\par\let\thelognumber\relax\let\thevolumenumber\relax
\let\thepapernumber\relax\let\thevolumeyear\relax\let\startpage\relax
\let\finishpage\relax\let\publishdate\relax\let\receiveddate\relax
\let\reviseddate\relax\let\accepteddate\relax\let\theasciititle\relax
\let\theasciiauthors\relax\let\theasciiaddress\relax
\let\theasciiabstract\relax

\let\theasciiemail\relax

\ifplaintex
\font\logobig=cmssbx10 scaled 3836
\font\logomed=cmssbx10 scaled 2557
\else
\font\logobig=cmssbx10 scaled 4200
\font\logomed=cmssbx10 scaled 2800
\fi

\long\def\makeagttitle{   
\count0=\startpage
\agt\hfill      
\hbox to 45truept{\vbox to 0pt{\vglue -13truept{\logomed A\kern -.37em{\logobig 
T}\kern -.38em G}\vss}\hss}
\break
{\small Volume \thevolumenumber\ (\thevolumeyear)
\startpage--\finishpage\nl
Published: \publishdate}

\vglue .25truein

{\parskip=0pt\leftskip 0pt plus
1fil\def\\{\par\smallskip}{\Large\bf\thetitle}\par\medskip} \vglue
0.05truein

{\parskip=0pt\leftskip 0pt plus 1fil\def\\{\par}{\sc\theauthors}
\par\medskip}%
 
\vglue 0.03truein 

{\small\leftskip 25truept\rightskip 25truept{\bf Abstract}\stdspace\theabstract

{\bf AMS Classification}\stdspace\theprimaryclass
\ifx\thesecondaryclass\relax\else; \thesecondaryclass\fi\par
{\bf Keywords}\stdspace \thekeywords\par}\vglue 7truept

}   

\ifplaintex
\font\phead=cmsl9 scaled 950
\font\pnum=cmbx10 scaled 913
\font\pfoot=cmsl9 scaled 950
\headline{\vbox to 0pt{\vskip -4.5mm\line{\small\phead\ifnum
\count0=\startpage ISSN 1472-2739 (on-line) 1472-2747 (printed)
\hfill {\pnum\folio}\else\ifodd\count0\def\\{ }%
\ifx\theshorttitle\relax\thetitle\else\theshorttitle\fi\hfill{\pnum\folio}
\else\def\\{ and }{\pnum\folio}\hfill\ifx\theshortauthors\relax\theauthors
\else\theshortauthors\fi\fi\fi}\vss}}
\footline{\vbox to 0pt{\vglue 0mm\line{\small\pfoot\ifnum\count0=\startpage
\copyright\ \gtp\hfill\else
\agt, Volume \thevolumenumber\ (\thevolumeyear)\hfill\fi}\vss}}
\else
\font=cmsl9 scaled 1050
\font\lnum=cmbx10 
\font=cmsl9 scaled 1050
\makeatletter
\def\@oddhead{{\small\lhead\ifnum\count0=\startpage ISSN 1472-2739 
(on-line) 1472-2747 (printed)\hfill {\lnum\number\count0}\else\ifodd\count0
\def\\{ }\ifx\theshorttitle\relax \thetitle \else\theshorttitle\fi\hfill
{\lnum\number\count0}\else\def\\{ and }{\lnum\number\count0}
\hfill\ifx\theshortauthors\relax 
\theauthors\else\theshortauthors\fi\fi\fi}}\def\@evenhead{\@oddhead}
\def\@oddfoot{\small\lfoot\ifnum\count0=\startpage\copyright\ \gtp\hfill\else
\agt, Volume \thevolumenumber\ (\thevolumeyear)\hfill\fi}
\def\@evenfoot{\@oddfoot}
\makeatother
\fi
\let\maketitlepage\makeagttitle
\let\makeshorttitle\maketitlepage
\let\maketitle\maketitlepage

\newwrite\gtoutfile
\long\gdef\makeheadfile{  
{\def\\{, }\def\s{ }
\immediate\openout\gtoutfile head.xxx
\immediate\write\gtoutfile{To: math@arxiv.org}
\immediate\write\gtoutfile{Subject: put OR rep NNNNN:ppppp}
\immediate\write\gtoutfile{--text follows this line--}
\immediate\write\gtoutfile{Proxy-for: \ifx\theasciiauthors\relax
\theauthors\else\theasciiauthors\fi\s<\ifx\theasciiemail\relax\theemail\else\theasciiemail\fi>}
\immediate\write\gtoutfile{\noexpand\\}
\immediate\write\gtoutfile{Authors: \ifx\theasciiauthors\relax
\theauthors\else\theasciiauthors\fi}
{\def\\{ }\immediate\write\gtoutfile{Title: \ifx\theasciititle\relax
\thetitle\else\theasciititle\fi}}
\immediate\write\gtoutfile{Subj-class: GT or SG, GR etc}
\immediate\write\gtoutfile{MSC-class: \theprimaryclass\ifx\thesecondaryclass\relax\else, \thesecondaryclass\fi}
\immediate\write\gtoutfile{Journal-ref: Algebr. Geom. Topol. \thevolumenumber\s
(\thevolumeyear) \startpage-\finishpage}
\immediate\write\gtoutfile{Comments: Published by Algebraic and
Geometric Topology at}
\immediate\write\gtoutfile{\s\s\s  http://www.maths.warwick.ac.uk/agt/AGTVol\thevolumenumber/agt-\thevolumenumber-\thepapernumber.abs.html}
\immediate\write\gtoutfile{\noexpand\\}
\immediate\write\gtoutfile{}
\ifx\theasciiabstract\relax
\immediate\write\gtoutfile{\theabstract}\else
\immediate\write\gtoutfile{\theasciiabstract}\fi
\immediate\write\gtoutfile{}
\immediate\write\gtoutfile{\noexpand\\}
\immediate\write\gtoutfile{}
\immediate\closeout\gtoutfile}}  

\def\maketitlepage{\makeagttitle\makeheadfile}
\let\makeshorttitle\maketitlepage
\let\maketitle\maketitlepage

\def\ifplaintex{\expandafter\ifx\csname documentclass\endcsname\relax}

\def\gtp{{\mathsurround=0pt\it $\cal G\mskip-2mu$eometry \&\ 
$\cal T\!\!$opology $\cal P\!$ublications}}  

\def\recd{{\small Received:\qua\receiveddate\ifx\reviseddate\relax
\else\qquad Revised:\qua\reviseddate\fi\par}} 

\def\lognumber#1{\def\thelognumber{#1}}
\def\volumenumber#1{\def\thevolumenumber{#1}}
\def\volumeyear#1{\def\thevolumeyear{#1}}
\def\papernumber#1{\def\thepapernumber{#1}}
\def\pagenumbers#1#2{\def\startpage{#1}\def\finishpage{#2}}
\def\published#1{\def\publishdate{#1}}

\def\received#1{\def\receiveddate{#1}}
\def\revised#1{\def\reviseddate{#1}}
\def\accepted#1{\def\accepteddate{#1}}

\def\asciiaddress#1{\def\theasciiaddress{#1}}

\long\def\asciiabstract#1{\long\def\theasciiabstract{#1}}

\let\\\par\let\thelognumber\relax\let\thevolumenumber\relax
\let\thepapernumber\relax\let\thevolumeyear\relax\let\startpage\relax
\let\finishpage\relax\let\publishdate\relax\let\receiveddate\relax
\let\reviseddate\relax\let\accepteddate\relax\let\theasciititle\relax
\let\theasciiauthors\relax\let\theasciiaddress\relax
\let\theasciiabstract\relax

\let\theasciiemail\relax

\ifplaintex
\font\logobig=cmssbx10 scaled 3836
\font\logomed=cmssbx10 scaled 2557
\else
\font\logobig=cmssbx10 scaled 4200
\font\logomed=cmssbx10 scaled 2800
\fi

\long\def\makeagttitle{   
\count0=\startpage
\agt\hfill      
\hbox to 45truept{\vbox to 0pt{\vglue -13truept{\logomed A\kern -.37em{\logobig 
T}\kern -.38em G}\vss}\hss}
\break
{\small Volume \thevolumenumber\ (\thevolumeyear)
\startpage--\finishpage\nl
Published: \publishdate}

\vglue .25truein

{\parskip=0pt\leftskip 0pt plus
1fil\def\\{\par\smallskip}{\Large\bf\thetitle}\par\medskip} \vglue
0.05truein

{\parskip=0pt\leftskip 0pt plus 1fil\def\\{\par}{\sc\theauthors}
\par\medskip}%
 
\vglue 0.03truein 

{\small\leftskip 25truept\rightskip 25truept{\bf Abstract}\stdspace\theabstract

{\bf AMS Classification}\stdspace\theprimaryclass
\ifx\thesecondaryclass\relax\else; \thesecondaryclass\fi\par
{\bf Keywords}\stdspace \thekeywords\par}\vglue 7truept

}   

\ifplaintex
\font\phead=cmsl9 scaled 950
\font\pnum=cmbx10 scaled 913
\font\pfoot=cmsl9 scaled 950
\headline{\vbox to 0pt{\vskip -4.5mm\line{\small\phead\ifnum
\count0=\startpage ISSN 1472-2739 (on-line) 1472-2747 (printed)
\hfill {\pnum\folio}\else\ifodd\count0\def\\{ }%
\ifx\theshorttitle\relax\thetitle\else\theshorttitle\fi\hfill{\pnum\folio}
\else\def\\{ and }{\pnum\folio}\hfill\ifx\theshortauthors\relax\theauthors
\else\theshortauthors\fi\fi\fi}\vss}}
\footline{\vbox to 0pt{\vglue 0mm\line{\small\pfoot\ifnum\count0=\startpage
\copyright\ \gtp\hfill\else
\agt, Volume \thevolumenumber\ (\thevolumeyear)\hfill\fi}\vss}}
\else
\font=cmsl9 scaled 1050
\font\lnum=cmbx10 
\font=cmsl9 scaled 1050
\makeatletter
\def\@oddhead{{\small\lhead\ifnum\count0=\startpage ISSN 1472-2739 
(on-line) 1472-2747 (printed)\hfill {\lnum\number\count0}\else\ifodd\count0
\def\\{ }\ifx\theshorttitle\relax \thetitle \else\theshorttitle\fi\hfill
{\lnum\number\count0}\else\def\\{ and }{\lnum\number\count0}
\hfill\ifx\theshortauthors\relax 
\theauthors\else\theshortauthors\fi\fi\fi}}\def\@evenhead{\@oddhead}
\def\@oddfoot{\small\lfoot\ifnum\count0=\startpage\copyright\ \gtp\hfill\else
\agt, Volume \thevolumenumber\ (\thevolumeyear)\hfill\fi}
\def\@evenfoot{\@oddfoot}
\makeatother
\fi
\let\maketitlepage\makeagttitle
\let\makeshorttitle\maketitlepage
\let\maketitle\maketitlepage

\newwrite\gtoutfile
\long\gdef\makeheadfile{  
{\def\\{, }\def\s{ }
\immediate\openout\gtoutfile head.xxx
\immediate\write\gtoutfile{To: math@arxiv.org}
\immediate\write\gtoutfile{Subject: put OR rep NNNNN:ppppp}
\immediate\write\gtoutfile{--text follows this line--}
\immediate\write\gtoutfile{Proxy-for: \ifx\theasciiauthors\relax
\theauthors\else\theasciiauthors\fi\s<\ifx\theasciiemail\relax\theemail\else\theasciiemail\fi>}
\immediate\write\gtoutfile{\noexpand\\}
\immediate\write\gtoutfile{Authors: \ifx\theasciiauthors\relax
\theauthors\else\theasciiauthors\fi}
{\def\\{ }\immediate\write\gtoutfile{Title: \ifx\theasciititle\relax
\thetitle\else\theasciititle\fi}}
\immediate\write\gtoutfile{Subj-class: GT or SG, GR etc}
\immediate\write\gtoutfile{MSC-class: \theprimaryclass\ifx\thesecondaryclass\relax\else, \thesecondaryclass\fi}
\immediate\write\gtoutfile{Journal-ref: Algebr. Geom. Topol. \thevolumenumber\s
(\thevolumeyear) \startpage-\finishpage}
\immediate\write\gtoutfile{Comments: Published by Algebraic and
Geometric Topology at}
\immediate\write\gtoutfile{\s\s\s  http://www.maths.warwick.ac.uk/agt/AGTVol\thevolumenumber/agt-\thevolumenumber-\thepapernumber.abs.html}
\immediate\write\gtoutfile{\noexpand\\}
\immediate\write\gtoutfile{}
\ifx\theasciiabstract\relax
\immediate\write\gtoutfile{\theabstract}\else
\immediate\write\gtoutfile{\theasciiabstract}\fi
\immediate\write\gtoutfile{}
\immediate\write\gtoutfile{\noexpand\\}
\immediate\write\gtoutfile{}
\immediate\closeout\gtoutfile}}  

\def\maketitlepage{\makeagttitle\makeheadfile}
\let\makeshorttitle\maketitlepage
\let\maketitle\maketitlepage

\def\ifplaintex{\expandafter\ifx\csname documentclass\endcsname\relax}

\def\gtp{{\mathsurround=0pt\it $\cal G\mskip-2mu$eometry \&\ 
$\cal T\!\!$opology $\cal P\!$ublications}}  

\def\recd{{\small Received:\qua\receiveddate\ifx\reviseddate\relax
\else\qquad Revised:\qua\reviseddate\fi\par}} 

\def\lognumber#1{\def\thelognumber{#1}}
\def\volumenumber#1{\def\thevolumenumber{#1}}
\def\volumeyear#1{\def\thevolumeyear{#1}}
\def\papernumber#1{\def\thepapernumber{#1}}
\def\pagenumbers#1#2{\def\startpage{#1}\def\finishpage{#2}}
\def\published#1{\def\publishdate{#1}}

\def\received#1{\def\receiveddate{#1}}
\def\revised#1{\def\reviseddate{#1}}
\def\accepted#1{\def\accepteddate{#1}}

\def\asciiaddress#1{\def\theasciiaddress{#1}}

\long\def\asciiabstract#1{\long\def\theasciiabstract{#1}}

\let\\\par\let\thelognumber\relax\let\thevolumenumber\relax
\let\thepapernumber\relax\let\thevolumeyear\relax\let\startpage\relax
\let\finishpage\relax\let\publishdate\relax\let\receiveddate\relax
\let\reviseddate\relax\let\accepteddate\relax\let\theasciititle\relax
\let\theasciiauthors\relax\let\theasciiaddress\relax
\let\theasciiabstract\relax

\let\theasciiemail\relax

\ifplaintex
\font\logobig=cmssbx10 scaled 3836
\font\logomed=cmssbx10 scaled 2557
\else
\font\logobig=cmssbx10 scaled 4200
\font\logomed=cmssbx10 scaled 2800
\fi

\long\def\makeagttitle{   
\count0=\startpage
\agt\hfill      
\hbox to 45truept{\vbox to 0pt{\vglue -13truept{\logomed A\kern -.37em{\logobig 
T}\kern -.38em G}\vss}\hss}
\break
{\small Volume \thevolumenumber\ (\thevolumeyear)
\startpage--\finishpage\nl
Published: \publishdate}

\vglue .25truein

{\parskip=0pt\leftskip 0pt plus
1fil\def\\{\par\smallskip}{\Large\bf\thetitle}\par\medskip} \vglue
0.05truein

{\parskip=0pt\leftskip 0pt plus 1fil\def\\{\par}{\sc\theauthors}
\par\medskip}%
 
\vglue 0.03truein 

{\small\leftskip 25truept\rightskip 25truept{\bf Abstract}\stdspace\theabstract

{\bf AMS Classification}\stdspace\theprimaryclass
\ifx\thesecondaryclass\relax\else; \thesecondaryclass\fi\par
{\bf Keywords}\stdspace \thekeywords\par}\vglue 7truept

}   

\ifplaintex
\font\phead=cmsl9 scaled 950
\font\pnum=cmbx10 scaled 913
\font\pfoot=cmsl9 scaled 950
\headline{\vbox to 0pt{\vskip -4.5mm\line{\small\phead\ifnum
\count0=\startpage ISSN 1472-2739 (on-line) 1472-2747 (printed)
\hfill {\pnum\folio}\else\ifodd\count0\def\\{ }%
\ifx\theshorttitle\relax\thetitle\else\theshorttitle\fi\hfill{\pnum\folio}
\else\def\\{ and }{\pnum\folio}\hfill\ifx\theshortauthors\relax\theauthors
\else\theshortauthors\fi\fi\fi}\vss}}
\footline{\vbox to 0pt{\vglue 0mm\line{\small\pfoot\ifnum\count0=\startpage
\copyright\ \gtp\hfill\else
\agt, Volume \thevolumenumber\ (\thevolumeyear)\hfill\fi}\vss}}
\else
\font=cmsl9 scaled 1050
\font\lnum=cmbx10 
\font=cmsl9 scaled 1050
\makeatletter
\def\@oddhead{{\small\lhead\ifnum\count0=\startpage ISSN 1472-2739 
(on-line) 1472-2747 (printed)\hfill {\lnum\number\count0}\else\ifodd\count0
\def\\{ }\ifx\theshorttitle\relax \thetitle \else\theshorttitle\fi\hfill
{\lnum\number\count0}\else\def\\{ and }{\lnum\number\count0}
\hfill\ifx\theshortauthors\relax 
\theauthors\else\theshortauthors\fi\fi\fi}}\def\@evenhead{\@oddhead}
\def\@oddfoot{\small\lfoot\ifnum\count0=\startpage\copyright\ \gtp\hfill\else
\agt, Volume \thevolumenumber\ (\thevolumeyear)\hfill\fi}
\def\@evenfoot{\@oddfoot}
\makeatother
\fi
\let\maketitlepage\makeagttitle
\let\makeshorttitle\maketitlepage
\let\maketitle\maketitlepage

\newwrite\gtoutfile
\long\gdef\makeheadfile{  
{\def\\{, }\def\s{ }
\immediate\openout\gtoutfile head.xxx
\immediate\write\gtoutfile{To: math@arxiv.org}
\immediate\write\gtoutfile{Subject: put OR rep NNNNN:ppppp}
\immediate\write\gtoutfile{--text follows this line--}
\immediate\write\gtoutfile{Proxy-for: \ifx\theasciiauthors\relax
\theauthors\else\theasciiauthors\fi\s<\ifx\theasciiemail\relax\theemail\else\theasciiemail\fi>}
\immediate\write\gtoutfile{\noexpand\\}
\immediate\write\gtoutfile{Authors: \ifx\theasciiauthors\relax
\theauthors\else\theasciiauthors\fi}
{\def\\{ }\immediate\write\gtoutfile{Title: \ifx\theasciititle\relax
\thetitle\else\theasciititle\fi}}
\immediate\write\gtoutfile{Subj-class: GT or SG, GR etc}
\immediate\write\gtoutfile{MSC-class: \theprimaryclass\ifx\thesecondaryclass\relax\else, \thesecondaryclass\fi}
\immediate\write\gtoutfile{Journal-ref: Algebr. Geom. Topol. \thevolumenumber\s
(\thevolumeyear) \startpage-\finishpage}
\immediate\write\gtoutfile{Comments: Published by Algebraic and
Geometric Topology at}
\immediate\write\gtoutfile{\s\s\s  http://www.maths.warwick.ac.uk/agt/AGTVol\thevolumenumber/agt-\thevolumenumber-\thepapernumber.abs.html}
\immediate\write\gtoutfile{\noexpand\\}
\immediate\write\gtoutfile{}
\ifx\theasciiabstract\relax
\immediate\write\gtoutfile{\theabstract}\else
\immediate\write\gtoutfile{\theasciiabstract}\fi
\immediate\write\gtoutfile{}
\immediate\write\gtoutfile{\noexpand\\}
\immediate\write\gtoutfile{}
\immediate\closeout\gtoutfile}}  

\def\maketitlepage{\makeagttitle\makeheadfile}
\let\makeshorttitle\maketitlepage
\let\maketitle\maketitlepage

\def\ifplaintex{\expandafter\ifx\csname documentclass\endcsname\relax}

\def\gtp{{\mathsurround=0pt\it $\cal G\mskip-2mu$eometry \&\ 
$\cal T\!\!$opology $\cal P\!$ublications}}  

\def\recd{{\small Received:\qua\receiveddate\ifx\reviseddate\relax
\else\qquad Revised:\qua\reviseddate\fi\par}} 

\def\lognumber#1{\def\thelognumber{#1}}
\def\volumenumber#1{\def\thevolumenumber{#1}}
\def\volumeyear#1{\def\thevolumeyear{#1}}
\def\papernumber#1{\def\thepapernumber{#1}}
\def\pagenumbers#1#2{\def\startpage{#1}\def\finishpage{#2}}
\def\published#1{\def\publishdate{#1}}

\def\received#1{\def\receiveddate{#1}}
\def\revised#1{\def\reviseddate{#1}}
\def\accepted#1{\def\accepteddate{#1}}

\def\asciiaddress#1{\def\theasciiaddress{#1}}

\long\def\asciiabstract#1{\long\def\theasciiabstract{#1}}

\let\\\par\let\thelognumber\relax\let\thevolumenumber\relax
\let\thepapernumber\relax\let\thevolumeyear\relax\let\startpage\relax
\let\finishpage\relax\let\publishdate\relax\let\receiveddate\relax
\let\reviseddate\relax\let\accepteddate\relax\let\theasciititle\relax
\let\theasciiauthors\relax\let\theasciiaddress\relax
\let\theasciiabstract\relax

\let\theasciiemail\relax

\ifplaintex
\font\logobig=cmssbx10 scaled 3836
\font\logomed=cmssbx10 scaled 2557
\else
\font\logobig=cmssbx10 scaled 4200
\font\logomed=cmssbx10 scaled 2800
\fi

\long\def\makeagttitle{   
\count0=\startpage
\agt\hfill      
\hbox to 45truept{\vbox to 0pt{\vglue -13truept{\logomed A\kern -.37em{\logobig 
T}\kern -.38em G}\vss}\hss}
\break
{\small Volume \thevolumenumber\ (\thevolumeyear)
\startpage--\finishpage\nl
Published: \publishdate}

\vglue .25truein

{\parskip=0pt\leftskip 0pt plus
1fil\def\\{\par\smallskip}{\Large\bf\thetitle}\par\medskip} \vglue
0.05truein

{\parskip=0pt\leftskip 0pt plus 1fil\def\\{\par}{\sc\theauthors}
\par\medskip}%
 
\vglue 0.03truein 

{\small\leftskip 25truept\rightskip 25truept{\bf Abstract}\stdspace\theabstract

{\bf AMS Classification}\stdspace\theprimaryclass
\ifx\thesecondaryclass\relax\else; \thesecondaryclass\fi\par
{\bf Keywords}\stdspace \thekeywords\par}\vglue 7truept

}   

\ifplaintex
\font\phead=cmsl9 scaled 950
\font\pnum=cmbx10 scaled 913
\font\pfoot=cmsl9 scaled 950
\headline{\vbox to 0pt{\vskip -4.5mm\line{\small\phead\ifnum
\count0=\startpage ISSN 1472-2739 (on-line) 1472-2747 (printed)
\hfill {\pnum\folio}\else\ifodd\count0\def\\{ }%
\ifx\theshorttitle\relax\thetitle\else\theshorttitle\fi\hfill{\pnum\folio}
\else\def\\{ and }{\pnum\folio}\hfill\ifx\theshortauthors\relax\theauthors
\else\theshortauthors\fi\fi\fi}\vss}}
\footline{\vbox to 0pt{\vglue 0mm\line{\small\pfoot\ifnum\count0=\startpage
\copyright\ \gtp\hfill\else
\agt, Volume \thevolumenumber\ (\thevolumeyear)\hfill\fi}\vss}}
\else
\font=cmsl9 scaled 1050
\font\lnum=cmbx10 
\font=cmsl9 scaled 1050
\makeatletter
\def\@oddhead{{\small\lhead\ifnum\count0=\startpage ISSN 1472-2739 
(on-line) 1472-2747 (printed)\hfill {\lnum\number\count0}\else\ifodd\count0
\def\\{ }\ifx\theshorttitle\relax \thetitle \else\theshorttitle\fi\hfill
{\lnum\number\count0}\else\def\\{ and }{\lnum\number\count0}
\hfill\ifx\theshortauthors\relax 
\theauthors\else\theshortauthors\fi\fi\fi}}\def\@evenhead{\@oddhead}
\def\@oddfoot{\small\lfoot\ifnum\count0=\startpage\copyright\ \gtp\hfill\else
\agt, Volume \thevolumenumber\ (\thevolumeyear)\hfill\fi}
\def\@evenfoot{\@oddfoot}
\makeatother
\fi
\let\maketitlepage\makeagttitle
\let\makeshorttitle\maketitlepage
\let\maketitle\maketitlepage

\newwrite\gtoutfile
\long\gdef\makeheadfile{  
{\def\\{, }\def\s{ }
\immediate\openout\gtoutfile head.xxx
\immediate\write\gtoutfile{To: math@arxiv.org}
\immediate\write\gtoutfile{Subject: put OR rep NNNNN:ppppp}
\immediate\write\gtoutfile{--text follows this line--}
\immediate\write\gtoutfile{Proxy-for: \ifx\theasciiauthors\relax
\theauthors\else\theasciiauthors\fi\s<\ifx\theasciiemail\relax\theemail\else\theasciiemail\fi>}
\immediate\write\gtoutfile{\noexpand\\}
\immediate\write\gtoutfile{Authors: \ifx\theasciiauthors\relax
\theauthors\else\theasciiauthors\fi}
{\def\\{ }\immediate\write\gtoutfile{Title: \ifx\theasciititle\relax
\thetitle\else\theasciititle\fi}}
\immediate\write\gtoutfile{Subj-class: GT or SG, GR etc}
\immediate\write\gtoutfile{MSC-class: \theprimaryclass\ifx\thesecondaryclass\relax\else, \thesecondaryclass\fi}
\immediate\write\gtoutfile{Journal-ref: Algebr. Geom. Topol. \thevolumenumber\s
(\thevolumeyear) \startpage-\finishpage}
\immediate\write\gtoutfile{Comments: Published by Algebraic and
Geometric Topology at}
\immediate\write\gtoutfile{\s\s\s  http://www.maths.warwick.ac.uk/agt/AGTVol\thevolumenumber/agt-\thevolumenumber-\thepapernumber.abs.html}
\immediate\write\gtoutfile{\noexpand\\}
\immediate\write\gtoutfile{}
\ifx\theasciiabstract\relax
\immediate\write\gtoutfile{\theabstract}\else
\immediate\write\gtoutfile{\theasciiabstract}\fi
\immediate\write\gtoutfile{}
\immediate\write\gtoutfile{\noexpand\\}
\immediate\write\gtoutfile{}
\immediate\closeout\gtoutfile}}  

\def\maketitlepage{\makeagttitle\makeheadfile}
\let\makeshorttitle\maketitlepage
\let\maketitle\maketitlepage

\lognumber{36}
\volumenumber{1}
\volumeyear{2001}
\papernumber{36}
\pagenumbers{719}{742}
\received{17 October 2000}
\revised{12 October 2001}
\accepted{26 Novemver 2001}
\published{1 December 2001}

\usepackage{amsmath,amssymb}

\usepackage{pb-diagram}

\numberwithin{equation}{section}

\newtheorem{theor}[equation]{Theorem}
\newtheorem{cor}[equation]{Corollary}
\newtheorem{lem}[equation]{Lemma}
\newtheorem{proposition}[equation]{Proposition}
\newtheorem{propriete}[equation]{Property}

\theoremstyle{definition}
\newtheorem{defin}[equation]{Definition}

\newtheorem{ex}[equation]{Example}

\newcommand{\TA}{\mathrm{T\!A}}

\def\build#1_#2^#3{\mathrel{\mathop{\kern0pt#1}\limits_{#2}^{#3}}}
\begin{document}
\title{On the cohomology algebra of a fiber}
\author{Luc Menichi}
\address{Universit\'e d'Angers, Facult\'e des Sciences\\2 
Boulevard Lavoisier, 49045 Angers, FRANCE}
\asciiaddress{Universite d'Angers, Faculte des Sciences\\2 
Boulevard Lavoisier, 49045 Angers, FRANCE}
\email{Luc.Menichi@univ-angers.fr}
\begin{abstract}
Let $f:E\rightarrow B$ be a fibration of fiber $F$.
Eilenberg and Moore have proved that there is a natural
isomorphism of vector spaces between $H^{*}(F;\mathbb{F}_p)$ and
$\mbox{Tor}^{C^{*}(B)}(C^{*}(E),\mathbb{F}_p)$.
Generalizing the rational case proved by Sullivan,
Anick \cite{AnickD:Hopah} proved that if $X$ is a finite $r$-connected
CW-complex of dimension $\leq rp$ then the algebra of singular cochains
$C^{*}(X;\mathbb{F}_p)$ can be replaced by a commutative differential graded
algebra $A(X)$ with the same cohomology.
Therefore if we suppose that $f:E\hookrightarrow B$ is an inclusion
of finite $r$-connected CW-complexes of dimension $\leq rp$,
we obtain an isomorphism of vector spaces between the algebra
 $H^{*}(F;\mathbb{F}_p)$ and
 $\mbox{Tor}^{A(B)}(A(E),\mathbb{F}_p)$
which has also a natural structure of algebra.
Extending the rational case proved by
Grivel-Thomas-Halperin~\cite{GrivelP:Fordss,Halperin:lecmm},
we prove that this isomorphism is in fact an isomorphism of
algebras. In particular, $H^{*}(F;\mathbb{F}_p)$ is a divided powers
algebra and $p^{th}$ powers vanish in the reduced cohomology
$\tilde{H}^{*}(F;\mathbb{F}_p)$.
\end{abstract}

\asciiabstract{Let f:E-->B be a fibration of fiber F.  Eilenberg and
Moore have proved that there is a natural isomorphism of vector spaces
between H^*(F;F_p) and Tor^{C^*(B)}(C^*(E),F_p).  Generalizing the
rational case proved by Sullivan, Anick [Hopf algebras up to homotopy,
J. Amer. Math. Soc. 2 (1989) 417--453] proved that if X is a finite
r-connected CW-complex of dimension < rp+1 then the algebra of
singular cochains C^*(X;F_p) can be replaced by a commutative
differential graded algebra A(X) with the same cohomology.  Therefore
if we suppose that f:E-->B is an inclusion of finite r-connected
CW-complexes of dimension < rp+1, we obtain an isomorphism of vector
spaces between the algebra H^*(F;F_p) and Tor^{A(B)}(A(E),F_p) which
has also a natural structure of algebra.  Extending the rational case
proved by Grivel-Thomas-Halperin [PP Grivel, Formes differentielles et
suites spectrales, Ann. Inst.  Fourier 29 (1979) 17--37] and [S
Halperin, Lectures on minimal models, Soc. Math. France 9-10 (1983)]
we prove that this isomorphism is in fact an isomorphism of
algebras. In particular, $H^*(F;F_p) is a divided powers algebra and
p-th powers vanish in the reduced cohomology \tilde{H}^*(F;F_p).}

\maketitle
\section{Introduction}
Let $f:E\twoheadrightarrow B$ be a fibration of fiber $F$ with simply connected base
$B$.
A major problem in Algebraic Topology is to compute the homotopy type of $F$.

In 1966, S. Eilenberg and J. Moore \cite{Eilenberg-Moore:homfcc} proved
that the cohomology of $F$ with coefficients in a field ${\Bbbk}$,
denoted $H^{*}(F;{\Bbbk})$, is entirely determined,
as graded vector spaces by the structure of $C^{*}(B;{\Bbbk})$-module
induced on $C^{*}(E;\Bbbk)$ through $f$.(Here $C^{*}(-;{\Bbbk})$
denotes the singular cochains.)
More precisely, they generalize the classical notion of derived functor ``Tor''
to the differential case and obtain a natural isomorphism of graded vector spaces
$$H^{*}(F)\cong\mbox{Tor}^{C^{*}(B)}(C^{*}(E),{\Bbbk}).$$

In the rational case, $C^{*}(X;\mathbb{Q})$ is equivalent to a
commutative cochain algebra
$A_{PL}(X)$~\cite{Sullivan:infct,Felix-Halperin-Thomas:ratht} which
carries the rational homotopy type of $X$. Moreover, the
Eilenberg-Moore isomorphism is induced by a quasi-isomorphism between
$A_{PL}(F)$ and a commutative cochain algebra $A$ constructed from
$A_{PL}(B)$ and $A_{PL}(E)$. In particular, the Eilenberg-Moore
isomorphism is an isomorphism of graded algebras.

In general, the Eilenberg-Moore isomorphism does not give the
multiplicative structure of $H^{*}(F;{\Bbbk})$. However the main
result of this paper asserts that for char ${\Bbbk}$
``sufficiently large '', the Eilenberg-Moore isomorphism is an
isomorphism of graded algebras with respect to a natural
multiplicative structure on Tor.

We now give the precise statement of our main result:

Over a field ${\Bbbk}$ of positive characteristic $p$,
Anick~\cite[Proposition 8.7(a)]{AnickD:Hopah}
proved that if $X$ is a finite $r$-connected CW-complex of
$\mbox{dimension } \leq rp$,
the algebra of singular cochains 
$C^{*}(X)$ is naturally linked to a commutative differential graded algebra $A(X)$
by morphisms of differential graded algebras inducing
isomorphisms in homology.
Therefore if we suppose that $f:E\hookrightarrow B$ is an inclusion
of finite  $r$-connected CW-complexes
of dimension $\leq rp$, we obtain the isomorphism of graded vector spaces
$$\mbox{Tor}^{C^{*}(B)}(C^{*}(E),{\Bbbk})\cong
\mbox{Tor}^{A(B)}(A(E),{\Bbbk}).$$
Thus the Eilenberg-Moore isomorphism becomes
$$H^{*}(F;{\Bbbk})\cong \mbox{Tor}^{A(B)}(A(E),{\Bbbk}).$$

Now since $A(B)$ and $A(E)$ are commutative,
$\mbox{Tor}^{A(B)}(A(E),{\Bbbk})$ has a natural structure of algebra.
We prove

\noindent{\bf Theorem A}\qua
{\sl Assume the characteristic of the field ${\Bbbk}$ is an odd prime $p$
and consider an inclusion $E\hookrightarrow B$ of finite
$r$-connected CW-complexes $(r\geq 1)$
of dimension $\leq rp$.
Then the Eilenberg-Moore isomorphism
$$H^{*}(F;{\Bbbk})\cong \mbox{Tor}^{A(B)}(A(E),{\Bbbk})$$
is an isomorphism of graded algebras.}

As corollary of this main result, we obtain

\noindent{\bf Theorem B}\qua{\rm (\ref{relative Wilkerson's conjecture})}\qua
{\sl Let $p$ be an odd prime.
Consider the homotypy fiber $F$ of an inclusion of finite $r$-connected CW-complexes
of dimension $\leq rp$.
Then the cohomology algebra $H^{*}(F;\mathbb{F}_p)$
is a divided powers algebra.
In particular, $p^{th}$ powers vanish in the reduced cohomology
$\tilde{H}^{*}(F;\mathbb{F}_p)$.}

In fact, Theorem A is a consequence of the model theorem
(Theorem \ref{passage aux models avec cofibration}) which establishes,
for any fibration $F\hookrightarrow E\twoheadrightarrow B$,
the existence of a coalgebra model
up to homotopy of the $\Omega E$-fibration
$$\Omega E\rightarrow \Omega B\rightarrow F.$$
Approaches concerning the general problem of computing the cohomology
algebra of a fiber, different from our model theorem,
are given in~\cite{Dupont-Hess:twitmf} and in~\cite{NdombolB:algcqc}.
We would like to express our gratitude to N. Dupont and to
S. Halperin for many useful discussions and suggestions which led to this
 work. We also thank the referee for significant simplifications.
This research was supported by the University of Lille (URA CNRS 751)
and by the University of Toronto (NSERC grants RGPIN 8047-98 and OGP000 7885).
\section{The two-sided bar construction}
We use the terminology of \cite{Felix-Halperin-Thomas:dgait}.
In particular, a quasi-isomorphism is denoted
$\buildrel{\simeq}\over\rightarrow$.
Let $A$ be an augmented differential graded algebra,
$M$ a right $A$-module, $N$ a left $A$-module.
Denote by $d_1$ the differential of the complex
$M\otimes T(s\overline{A})\otimes N$ obtained by tensorization,
and denote by $s\overline{A}$ the suspension of the augmentation ideal
$\overline{A}$, $(s\overline{A})_i=\overline{A}_{i-1}$.
Let $\vert x\vert$ be the degree of an element $x$ in any graded object.
We denote the tensor product of the elements $m\in M$, $sa_1\in s\overline{A}$,
\ldots , $sa_k\in s\overline{A}$ and $n\in N$ by $m[sa_1|\cdots|sa_k]n$.
Let $d_2$ be the differential on the graded vector space
$M\otimes T(s\overline{A})\otimes N$ defined by:
\begin{eqnarray*}
d_2m[sa_1|\cdots|sa_k]n&=&(-1)^{\vert m\vert } ma_1[sa_2|\cdots|sa_k]n\\
&&+\sum _{i=1}^{k-1} (-1)^{\varepsilon_i}{m[sa_1|\cdots|sa_ia_{i+1}|\cdots|sa_k]n}\\
&&-(-1)^{\varepsilon_{k-1}} m[sa_1|\cdots|sa_{k-1}]a_kn;
\end{eqnarray*}
Here $\varepsilon_i=\vert m\vert +\vert sa_1\vert+\cdots +\vert sa_i\vert$.

The {\it bar construction of $A$ with coefficients in $M$ and $N$}, denoted
$B(M;A;N)$, is the complex $(M\otimes T(s\overline{A})\otimes N,d_1+d_2)$.
We use mainly $B(M;A)=B(M;A;{\Bbbk})$. The {\it reduced bar construction of $A$},
denoted $B(A)$, is $B({\Bbbk};A)$.

Let $B$ be another augmented differential graded algebra, $P$ a right $B$-module and
$Q$ a left $B$-module.
Then we have the natural Alexander-Whitney morphism of complexes
(\cite[X.7.2]{MacLane:hom} or \cite[XI.6(3) computation of the
$\vee$ product]{Cartan-Eilenberg:homalg})
$$AW:B(M\otimes P;A \otimes B;N\otimes Q)\rightarrow B(M;A;N)\otimes B(P;B;Q)$$
where the image of a typical element
$m\otimes p [s(a_1 \otimes b_1)|\cdots |s(a_k \otimes b_k)]n\otimes q$ is\\
\centerline{$\displaystyle\sum _{i=0}^{k}(-1)^{\zeta_i}m[sa_1|\cdots |sa_i]a_{i+1}\cdots a_kn
\otimes pb_1\cdots b_i[sb_{i+1}|\cdots |sb_k]q.$}
\begin{eqnarray*}
\mbox{Here}\quad\zeta_i&=&\sum_{j=1}^{k}\left (\vert p\vert+\sum_{l=1}^{j-1}\vert b_l\vert\right )
\vert a_j\vert+\left (\vert p\vert+\sum_{j=1}^{k}\vert b_j\vert\right )\vert n\vert\\
&&+\sum_{j=i+1}^{k}(j-i)\vert a_j\vert+(k-i)\vert n\vert+\vert i\vert\vert p\vert
+\sum_{j=1}^{i-1}(i-j)\vert b_j\vert.
\end{eqnarray*}
\begin{propriete}\label{structure de coalgebre sur la bar}
If there exist a morphism of augmented algebras
$\Delta_A :A\rightarrow A\otimes A$
and morphisms of $A$-modules
$\Delta_M:M\rightarrow M\otimes M$,
$\varepsilon_M:M\twoheadrightarrow{\Bbbk}$,
$\Delta_N:N\rightarrow N\otimes N$,
$\varepsilon_N:N\twoheadrightarrow{\Bbbk}$
then
$B(\varepsilon_M;\varepsilon_A;\varepsilon_N):B(M;A;N)\twoheadrightarrow
B({\Bbbk};{\Bbbk};{\Bbbk})={\Bbbk}$ is an augmentation for $B(M;A;N)$
and the composite
\begin{align*}
B(M;A;N)\xrightarrow {B(\Delta_M;\Delta_A;\Delta_N)}
{B(M\otimes M;A \otimes A;N\otimes N)}&\\
\buildrel {AW}\over \longrightarrow B(M;A;N)&\otimes B(M;A;N)
\end{align*}
is a morphism of complexes.
In particular, if $A$ is a differential graded Hopf algebra and if $M$ and $N$ are
$A$-coalgebras
then $B(M;A;N)$ is a differential graded coalgebra.
This coalgebra structure on $B(M;A;N)$ is functorial.
\end{propriete}
\begin{propriete}\label{bar construction d'un semi-libre}
Moreover, if $M$ is $A$-semifree (in the sense of \cite[\S 2]{Felix-Halperin-Thomas:dgait})
then $B(M;A;N)\buildrel {\simeq}\over\rightarrow M\otimes_A N$
is a quasi-isomorphism of coalgebras.
\end{propriete}
\begin{theor}{\rm(\cite[5.1]{Felix-Halperin-Thomas:dgait} or
\cite[Theorem 3.9 and Corollary 3.10]{Gugenheim-May:theadt})}\label{base}\qua
Let $p:E\twoheadrightarrow B$ be a right $G$-fibration with $B$ path connected.
Then there is a natural quasi-isomorphism of coalgebras
$$B(C_*(E);C_*(G))\buildrel {\simeq}\over\rightarrow C_*(B).$$
\end{theor}
\begin{cor}\label{fibre}
Let $f:E\rightarrow B$ be a continuous pointed map with $E$ and
$B$ path connected. If its homotopy fiber $F$ is path connected,
then there is a chain coalgebra $G(f)$ equipped with two natural
isomorphisms of chain coalgebras
$$C_*(F)\buildrel {\simeq}\over\leftarrow G(f)
\buildrel {\simeq}\over\rightarrow
B(C_*(\Omega B);C_*(\Omega E)).$$
\end{cor}
This Corollary proves that the cohomology algebra $H^{*}(F)$
is determined by the Hopf algebra morphism
$C_*(\Omega f):C_*(\Omega E)\rightarrow C_*(\Omega B)$.
This is the starting observation of our paper.
In the next section,
we extend Property \ref{structure de coalgebre sur la bar} to
Hopf algebras and coalgebras up to homotopy: i.e. we do not require strict
coassociativity of the diagonals.
\section{Hopf algebras and coalgebras up to homotopy}\label{Bar
  construction and homotopies}
Let $f$, $g:A\rightarrow B$ be two morphisms of augmented differential
graded algebras.
A linear map $h:A\rightarrow B$ of (lower)
degree $+1$ is a {\it homotopy of
algebras} from $f$ to $g$ denoted $h:f\thickapprox_a g$
if $\varepsilon_B h=0$, $hd+dh=f-g$ and
$h(xy)=h(x)g(y)+(-1)^{\vert x\vert}f(x)h(y)$ for $x,y\in A$. 
The symbol $\thickapprox $ will be  reserved
to the usual notion of chain homotopy.

A {\it (cocommutative) coalgebra up to homotopy} is a complex $C$ equipped with
a morphism of complexes $\Delta :C\rightarrow C\otimes C$
and a morphism $\varepsilon:C\rightarrow {\Bbbk}$
such that $(\varepsilon\otimes 1)\circ\Delta =1= (1\otimes
\varepsilon)\circ\Delta$ (strict counitary),
$(\Delta\otimes 1)\circ\Delta \thickapprox (1\otimes
\Delta)\circ\Delta$ (homotopy coassociativity) and
$\tau\circ\Delta\thickapprox\Delta$ (homotopy cocommutativity).
Here $\tau(x\otimes y)=(-1)^{\vert x\vert\vert y\vert}y\otimes x$.
Let $C$ and $C'$ be two coalgebras up to homotopy.
A morphism of complexes $f:C\rightarrow C'$ is a {\it
  morphism of coalgebras up to homotopy} if $\Delta f\thickapprox
(f\otimes f)\Delta$ and $\varepsilon\circ f=\varepsilon$.

A {\it (cocommutative) Hopf algebra up to homotopy}
is a differential graded algebra
$K$ equipped with two morphisms of algebras $\Delta :K\rightarrow
K\otimes K$ and $\varepsilon: K\rightarrow {\Bbbk}$ such
that $(\varepsilon\otimes 1)\circ\Delta =1= (1\otimes \varepsilon)\circ\Delta$,
$(\Delta\otimes 1)\circ\Delta \thickapprox_a (1\otimes \Delta)\circ\Delta$ and
$\tau\circ\Delta\thickapprox_a\Delta$.
Let $K$, $K'$ be two Hopf algebras up to homotopy.
A morphism of augmented differential graded algebras
$f:K\rightarrow K'$ is a {\it a
  morphism of Hopf algebras up to homotopy} if
$\Delta f\thickapprox_a (f\otimes f)\Delta$.

\begin{lem}\label{naturalite des homotopies}
Suppose $\varphi\thickapprox_a\varphi':A\rightarrow A'$ and
$\Psi\thickapprox_a\Psi ':M\rightarrow M'$ via algebraic homotopies
$h$ and $h'$, with
$\varphi$, $\varphi '$, $\Psi$, $\Psi '$ morphisms of
augmented chain algebras.
Let $f:A\rightarrow M$ and $g:A'\rightarrow M'$ be two morphisms
of augmented chain algebras such that $\Psi\circ f=g\circ\varphi$
and $\Psi'\circ f=g\circ\varphi'$.
We summarize this situation by the ``diagram"
\[
\begin{diagram}
\node{A} \arrow{e,t}{h:\varphi\thickapprox_a\varphi'}
\arrow{s,l}{f} \node{A'} \arrow{s,r}{g} \\
\node{M} \arrow{e,b}{h':\Psi\thickapprox_a\Psi '} \node{M'}
\end{diagram}
\]
 If $h'\circ f=g\circ h$ ({\it naturality of the homotopies})
 then the morphisms of augmented chain complexes $B(\Psi;\varphi)$ and
$B(\Psi ';\varphi ')$ are chain homotopic.
\end{lem}
\begin{proof}The explicit chain homotopy $\Theta$ between
$B(\Psi;\varphi)$ and $B(\Psi ';\varphi ')$ is given
by
\begin{multline*}
\Theta (m[sa_1|\dots|sa_k])=
h'(m)[s\varphi'(a_1)|\dots|s\varphi'(a_k)]\\
 -\sum_{i=1}^{k}(-1)^{\varepsilon_{i-1}}\Psi(m)
[s\varphi(a_1)|\dots|s\varphi(a_{i-1})|sh(a_i)|s\varphi'(a_{i+1})|\dots
|s\varphi'(a_k)]
\end{multline*}
Recall that $\varepsilon_{i-1}=\vert m\vert +\vert sa_1\vert+\cdots +\vert sa_{i-1}\vert$.
Since $\Theta$ is just the chain homotopy obtained by tensorization,
$$d_1\Theta+\Theta d_1=\Psi\otimes T(s\varphi)-\Psi'\otimes T(s\varphi').$$
It remains to check that $d_2\Theta+\Theta d_2=0$.\end{proof}

From Lemma \ref{naturalite des homotopies}, one deduces:
\begin{lem}\label{Hopf algebras a homotopie pres}
\begin{itemize}
\item[\rm(i)] Let $K$ (respectively $C$) be a Hopf algebra up to homotopy,
coassociative up to a homotopy $h_{assocK}$ (respectively $h_{assocC}$):
$(\Delta\otimes 1)\circ\Delta \thickapprox_a (1\otimes
\Delta)\circ\Delta$
and cocommutative up a homotopy $h_{comK}$ (respectively $h_{comC}$):
$\tau\circ\Delta\thickapprox_a\Delta$.
Let $f:K\rightarrow C$ be a morphism of augmented algebras
such that $\Delta_Cf=(f\otimes f)\Delta_K$,
$h_{assocC}f=(f\otimes f\otimes f)h_{assocK}$ and
$h_{comC}f=(f\otimes f)h_{comK}$
($f$ commutes with the diagonals and the homotopies of
coassociativity and cocommutativity).
Then $B(C;K)$ with the diagonal
$$B(C;K)\xrightarrow{B(\Delta_C;\Delta_K)}B(C\otimes C;K\otimes K)
\buildrel {AW}\over\longrightarrow  B(C;K)\otimes B(C;K)$$
is a (cocommutative) coalgebra up to homotopy.
\item[\rm(ii)]Suppose given the following cube of augmented chain algebras
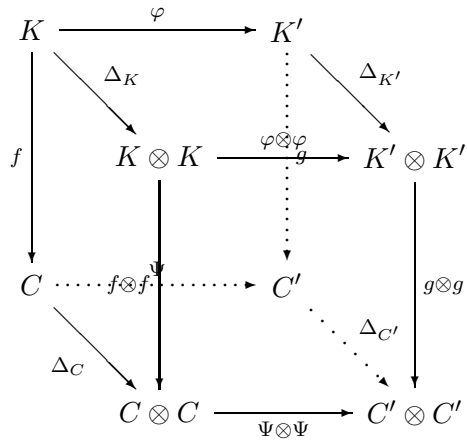
\[
\begin{diagram}[H]
       \node{K}   \arrow[2]{e,t}{\varphi}   \arrow[1]{se,t}{\Delta_K}   \arrow[2]{s,l}{f}     \node[2]{K'}  
                                                            \arrow[1]{se,t}{\Delta_{K'}} \arrow[2]{s,r,..}{g}\\
        \node{}   \node[1]{K\otimes K}             \arrow[2]{e,t}{\varphi\otimes\varphi}         \arrow[2]{s,l}{f\otimes f}             \node[2]{K'\otimes K'} 
                                                           \arrow[2]{s,r}{g\otimes g}\\
      \node{C}\arrow[2]{e,t,..}{\Psi}           \arrow[1]{se,b}{\Delta_C}   \node[2]{C'}\arrow[1]{se,t,..}{\Delta_{C'}}\\
         \node{}   \node[1]{C\otimes C}                       \arrow[2]{e,b}{\Psi\otimes\Psi}             \node[2]{C'\otimes C'}
\end{diagram}
\]
where all the faces commute exactly except the top and the
bottom ones.
Suppose that the top face commutes up to a homotopy
$h_{top}:(\varphi\otimes\varphi)\Delta_K\thickapprox_a\Delta_{K'}\varphi$
and the bottom face commutes up to a homotopy
$h_{bottom}:(\Psi\otimes\Psi)\Delta_C\thickapprox_a\Delta_{C'}\Psi$
such that $h_{bottom}f=(g\otimes g)h_{top}$.
Then the morphism of augmented chain complexes
$B(\Psi;\varphi):B(C;K)\rightarrow B(C';K')$ commutes with the
diagonals up to chain homotopy.
\end{itemize}
\end{lem}
\section{The model Theorem}\label{HAH structure on free models}
Let $X$ be a graded vector space. We denote a free chain
algebra $(TX,\partial)$
simply by $TX$ except when the differential $\partial$ can be specified.
In particular, a free chain algebra with zero differential is still denoted
by $(TX,0)$.

Let $f:E\rightarrow B$ be a map between path connected pointed topological spaces
with a path connected homotopy fiber $F$.
Then there is a commutative diagram of augmented chain algebras as
follows:
\begin{equation}\label{diagram du model theorem}
\begin{diagram}
\node{TX}\arrow{e,tb}{\simeq}{m_X}\arrow{s,l}{m(f)}
\node{C_*(\Omega E)}\arrow{s,r}{C_*(\Omega f)}\\
\node{TY}\arrow{e,tb}{\simeq}{m_Y}\node{C_*(\Omega B)}
\end{diagram}
\end{equation}
where
$TX$, $TY$ are free chain algebras, $m_X$, $m_Y$ are quasi-isomorphisms
and $m(f):TX \rightarrowtail TY$ is a free extension
(in the sense of \cite[\S 3]{Felix-Halperin-Thomas:dgait}).
\begin{theor}\label{passage aux models avec cofibration}
With the above:\begin{enumerate}
\item $TX$ (respectively $TY$) can be endowed with an structure of
Hopf algebras up to homotopy
such that $m_X$ (respectively $m_Y$) commutes with the diagonals up to
a homotopy $h_{X}$ (respectively $h_{Y}$) and such that
the diagonal of $TY$ extends the diagonal of $TX$,
the homotopy of coassociativity of $TY$ extends
the homotopy of coassociativity of $TY$,
the homotopy of cocommutativity of $TY$ extends
the homotopy of cocommutativity of $TX$
and $h_{Y}$ extends
$(C_*(\Omega f)\otimes C_*(\Omega f))h_{X}$.
\item $B(m_Y;m_X):B(TY;TX)\buildrel {\simeq}\over\rightarrow
B(C_*(\Omega B);C_*(\Omega E))$ is a
morphism of coalgebras up to homotopy.
\item The homology of the coalgebra up to homotopy $TY\otimes_{TX} {\Bbbk}$ is isomorphic to $H_*(F)$
as coalgebras.
\end{enumerate}
\end{theor}
It is easy to see that the isomorphism of graded coalgebras between $H_*(TY\otimes_{TX} {\Bbbk})$
and $H_*(F)$ fits into the commutative diagram of graded coalgebras:
\[
\begin{diagram}
\node{H_*(TY)}\arrow[2]{e,tb}{\cong}{H_*(m_Y)}\arrow{s,l}{H_*(q)}
\node[2]{H_*(\Omega B)}\arrow{s,r}{H_*(\partial)}\\
\node{H_*(TY\otimes_{TX} {\Bbbk})}\arrow[2]{e,t}{\cong}
\node[2]{H_*(F)}
\end{diagram}
\]
where $\partial:\Omega B\hookrightarrow F$ is the inclusion
$\Omega B\times *\subset PB\times_B E$
and $q:TY\twoheadrightarrow TY\otimes_{TX} {\Bbbk}$ is the quotient map.

The exact commutativity of the
diagram \ref{diagram du model theorem}
is not important. If the diagram commutes only up to homotopy,
since $m(f)$ is a free extension,
using the lifting lemma \cite[3.6]{Felix-Halperin-Thomas:dgait},
we can replace $m_Y$ by another $m_Y$ which is homotopic to it, so that
now the diagram strictly commutes.
On the contrary, it is important that $m(f)$ is a free extension.
We will show it in Section
\ref{fibre d'une application suspendue et contreexemple}.
Indeed, the general idea for the proof of part 1.
is to keep control of the homotopies
using the homotopy extension property of cofibrations.

\begin{proof}[Proof of Theorem \ref{passage aux models avec cofibration}]
1. By the lifting lemma (\cite[I.7 and II.1.11=II.2.11a)]{Baues:algh}
or \cite[3.6]{Felix-Halperin-Thomas:dgait}),
we obtain a diagonal $\Delta_{TX}$ for $TX$
such that the following diagram of augmented augmented chain algebras commutes up to
a homotopy $h_{X}$:
\[
\begin{diagram}
\node{TX}\arrow[2]{e,tb}{\simeq}{m_X}\arrow{s,l,..}{\Delta_{TX}}
\node[2]{C_*(\Omega E)}\arrow{s,r}{\Delta_{C_*(\Omega E)}}\\
\node{TX\otimes TX}\arrow[2]{e,tb}{\simeq}{m_X\otimes m_X}
\node[2]{C_*(\Omega E)\otimes C_*(\Omega E)}
\end{diagram}
\]
Moreover, since
$C_*(\Omega E)$ is a differential graded
Hopf algebra which is cocommutative up to
homotopy~\cite[Proposition 7.1]{AnickD:Hopah},
by the unicity of the lifting (\cite[II.1.11c)]{Baues:algh}
or \cite[3.7]{Felix-Halperin-Thomas:dgait}), $\Delta_{TX}$ is
counitary, coassociative and cocommutative, all up to
homotopy.
The diagonal $\Delta_{TX}$ can be chosen to be stricly
counitary~\cite[Lemma 5.4]{AnickD:Hopah}.
So $TX$ is an Hopf algebra up to homotopy.

By the naturality of the lifting lemma with respect to the inclusion
$m(f):TX \rightarrowtail TY$ \cite[3.6]{Felix-Halperin-Thomas:dgait},
we may put a diagonal on $TY$, $\Delta_{TY}$
extending the diagonal on $TX$ and
there exists a homotopy $h_{Y}$ between $(m_Y\otimes m_Y)\Delta_{TY}$
and $\Delta_{C_*(\Omega B)}m_Y$ extending $C_*(\Omega f)^{\otimes
  2}h_{X}$.
Again the diagonal on $TY$ can be chosen to be counitary and so $TY$
is also a Hopf algebra up to homotopy.

We give now a detailed proof that $\Delta_{TX}$ is cocommutative up to a homotopy
$h_{comX}$ and that  $\Delta_{TY}$ is cocommutative up to a homotopy $h_{comY}$
extending  $h_{comX}$.
Since the diagonal on $C_*(\Omega E)$ is cocommutative up to a homotopy
$h_{comE}$, 
by the unicity of the lifting (\cite[II.1.11c)]{Baues:algh}
or~\cite[3.7]{Felix-Halperin-Thomas:dgait}), $\Delta_{TX}$ is
cocommutative up to a homotopy $h_{comX}$.
More precisely,
$h_{comX}$ is such that in the diagram 
\[
\begin{diagram}[HHHHH]
\node{TX\amalg TX}\arrow[2]{e,t}{(\Delta_{TX},\tau\Delta_{TX})}\arrow{s,l}{(i_0,i_1)}
\node[2]{TX^{\otimes 2}}\arrow{s,lr}{m_X}{\simeq}\\
\node{ITX}\arrow{ene,t,..}{h_{comX}}
\arrow[2]{e,b}{h_X-h_{comE}\circ m_X-\tau h_X}
\node[2]{C_*(\Omega E)^{\otimes 2}}
\end{diagram}
\]
where $ITX$ is the Baues-Lemaire
cylinder (\cite[3.5]{Felix-Halperin-Thomas:dgait} or
\cite[I.7.12]{Baues:algh}),
the upper triangle commutes~\cite[II.1.11a)]{Baues:algh} and the lower
triangle commutes up to a homotopy relative to
$TX\amalg TX$~\cite[II.1.11b)]{Baues:algh}.
Now, since the homotopy of cocommutativity of $C_*(\Omega B)$ is
natural~\cite[(23)]{AnickD:Hopah} and the sums and negatives of homotopies are canonically defined~\cite[II.17.3]{Baues:algh},
the homotopy $h_X-h_{comE}\circ m_X-\tau h_X$ is extended by
$h_Y-h_{comB}\circ m_Y-\tau h_Y$.
Therefore, by push out, we obtain a morphism 
$ITX\cup_{TX\amalg TX}(TY\amalg TY)\rightarrow TY^{\otimes 2}$
extending $m(f)^{\otimes 2}\circ h_{comX}$, $\Delta_{TY}$
and $\tau\Delta_{TY}$.
The following square of unbroken arrows commutes up to homotopy:
\[
\begin{diagram}[HHHH]
\node{ITX\cup_{TX\amalg TX}(TY\amalg TY)}\arrow{s,l}{(I(m(f)),i_o,i_1)}
\arrow[2]{e}\node[2]{TY^{\otimes 2}}\arrow{s,lr}{m_Y}{\simeq}\\
\node{ITY}\arrow[2]{e,b}{h_Y-h_{comB}\circ m_Y-\tau h_Y}
\arrow{ene,b,..}{h_{comY}}
\node[2]{C_*(\Omega B)^{\otimes 2}}
\end{diagram}
\]
Using again the naturality of the
lifting lemma~\cite[3.6]{Felix-Halperin-Thomas:dgait}, we obtain
the homotopy of cocommutativity of $TY$, $h_{comY}$.
A similar proof shows that the homotopy of coassociativity on
$TY$ can be chosen to extend the homotopy of
coassociativity
on $TX$.
So finally, the whole structure (homotopies included) of Hopf
algebra up to homotopy on $TY$ extends the structure on $TX$
(Compare with the proof of Theorem 8.5(g)\cite{AnickD:Hopah}).

2. Now Lemma \ref{Hopf algebras a homotopie pres} says exactly that
part 1. implies part 2.

3. Since $TX\rightarrowtail TY$ is a semi-free extension of
$TX$-modules
(in the sense of \cite[\S 2]{Felix-Halperin-Thomas:dgait})
and by Property \ref{bar construction d'un semi-libre},
the quasi-isomorphism of augmented chain complexes
$$B(TY;TX)\buildrel{\simeq}\over\rightarrow TY\otimes_{TX} {\Bbbk}$$
commutes exactly with the diagonals.

Since $B(TY;TX)\buildrel{\simeq}\over\rightarrow TY\otimes_{TX}
{\Bbbk}$
is a diagonal preserving chain homotopy equivalence,
the diagonal in $TY\otimes_{TX} {\Bbbk}$ is homotopy coassociative
and homotopy cocommutative.
By Corollary \ref{fibre}, $C_*(F)$ is weakly equivalent to
$B(C_*\Omega B;C_*\Omega E)$ as coalgebras
So part 2. implies that the coalgebra $H_*(TY\otimes_{TX} {\Bbbk})$ is isomorphic to $H_*(F)$.\end{proof}
\section{The fiber of a suspended map}\label{fibre d'une application suspendue et contreexemple}
Let $C$ be a coaugmented differential graded coalgebra.
Consider the tensor algebra on $\overline{C}$, $T\overline{C}$,
equipped with the differential obtained by tensorisation.
The composite
$C\buildrel{\Delta_C}\over\rightarrow C\otimes C\hookrightarrow
T\overline{C}\otimes T\overline{C}$ extends to an unique morphism
of augmented differential graded algebras
$$\Delta_{T\overline{C}}:T\overline{C}\rightarrow T\overline{C}\otimes
T\overline{C}.$$
The tensor algebra $T\overline{C}$ equipped with this structure of
differential graded Hopf algebra, is called
the {\it Hopf algebra obtained by tensorization of the
coalgebra $C$} and is denoted $\TA\overline{C}$ in this section.
\begin{lem}\label{consequence de Husemoller}
Let $X$ be a path connected space. Then there is a natural
quasi-isomorphism of Hopf algebras
$\TA\overline{C_*(X)}\buildrel{\simeq}\over\rightarrow C_*(\Omega\Sigma X)$.
\end{lem}
\begin{proof}
The adjunction map $ad$ induces a morphism of coaugmented coalgebras
$C_*(ad):C_*(X)\rightarrow C_*(\Omega\Sigma X)$.
By universal property of $\TA\overline{C_*(X)}$, $C_*(ad)$
extends to a natural morphism of Hopf algebras.
By the Bott-Samelson Theorem~\cite[appendix 2 Theorem 1.4]{Husemoller:fibb},
it is a quasi-isomorphism, since the functors $H$ and $T$
commute.\end{proof}

\begin{theor}\label{fibre d'une suspension}
Let $f:E\rightarrow B$ be a continuous map between path connected
spaces. Let $F$ be the homotopy fiber of ${\Sigma f}$.
Then $C_*(F)$ is naturally weakly equivalent as coalgebras to
$B(\TA\overline{C_*(B)};\TA\overline{C_*(E)})$. In particular, the
algebra $H^{*}(F)$ depends only of the morphism of coalgebras $C_*(f)$.
\end{theor}
\begin{proof} This is a direct consequence of Lemma~\ref{consequence de
  Husemoller}, Corollary~\ref{fibre} and
Property~\ref{structure de coalgebre sur la bar}.\end{proof}

Consider the homotopy commutative diagram of chain algebras.
\begin{equation}\label{diagram de suspension}
\begin{diagram}
\node{(TH_+(E),0)}\arrow[2]{e,tb}{\simeq}{m_X}\arrow{s,l}{TH_+(f)}
\node[2]{\TA\overline{C_*(E)}}\arrow{s,r}{\TA\overline{C_*(f)}}\\
\node{(TH_+(B),0)}\arrow[2]{e,tb}{\simeq}{m_Y}
\node[2]{\TA\overline{C_*(B)}}
\end{diagram}
\end{equation}
where $m_X$ and $m_Y$ induce the identity in homology.

When $H_*(f)$ is injective then $TH_+(f)$ is a free extension
and we may choose $m_Y$ so that the diagram~(\ref{diagram de suspension}) commutes.
Thus Theorem~\ref{passage aux models avec cofibration} applies.
The structures of Hopf algebra up to homotopy on $TH_+(E)$ and
$TH_+(B)$ given by part 1. of Theorem \ref{passage aux models
  avec cofibration} are the structures of Hopf algebra
obtained by tensorization of the coalgebras $H_*(E)$ and $H_*(B)$.
Part 3. of Theorem \ref{passage aux models
  avec cofibration} claims that we have the isomorphism of
graded coalgebras
\begin{equation}\label{fibre d'une injection suspendue}
\TA H_+(B)\otimes_{\TA H_+(E)}{\Bbbk}\cong H_{*}(F).
\end{equation}

If $H_*(f)$ is not injective, this is
not true in general: the algebra $H^{*}(F)$ does not depend only on
the
morphism of coalgebras $H_*(f)$.
Indeed, over $\mathbb{F}_p$,
take $f$ to be a map from
$S^{2p-1}$ to $\mathbb{CP} ^{p-1}$.
Whatever is the map chosen, $H_*(f):H_*(S^{2p-1})\rightarrow
H_*(\mathbb{CP} ^{p-1})$ is null.
Let $y_2$ be a generator of $H^{2}(F)$.
If $f$ is the Hopf map, there is a map
$\psi :\mathbb{CP} ^{p}\rightarrow F$ such that 
the following diagram commutes
\begin{equation*}
\begin{diagram}
\node{S^{2p-1}}\arrow[2]{e,t}{f}\arrow{s,t}{ad}
\node[2]{\mathbb{CP} ^{p-1}}\arrow[2]{e}\arrow{s,t}{ad}
\node[2]{\mathbb{CP} ^{p}}\arrow{s,b,..}{\psi}\\
\node{\Omega\Sigma S^{2p-1}}\arrow[2]{e,b}{\Omega\Sigma f}
\node[2]{\Omega\Sigma \mathbb{CP} ^{p-1}}\arrow[2]{e,b}{\partial}
\node[2]{F}
\end{diagram}
\end{equation*}
Since 
$H^{2}(\psi )$ is an isomorphism, $y_2^{p}
\not= 0$.
On the contrary, if $f$ is the constant map then
$F\thickapprox \Omega\Sigma \mathbb{CP} ^{p-1}\times\Sigma S^{2p-1}$
and so $y_2^{p}=0$.

Of course,
the isomorphism of coalgebras (\ref{fibre d'une injection suspendue})
can be proved more easily with the Eilenberg-Moore
spectral sequence applied to the $\Omega\Sigma E$-fibration
$$\Omega\Sigma E\buildrel{\Omega\Sigma f}\over\longrightarrow
\Omega\Sigma B\rightarrow F.$$
\section{Proof of Theorem A}
We recall first the natural structure of algebra on the torsion product of
commutative algebras.
Let $f:A\rightarrow M$, $g:A\rightarrow N$ be two morphisms of
commutative differential graded algebras.
The composite
$$\mbox{Tor}^{A}(M,N)\otimes\mbox{Tor}^{A}(M,N)\buildrel{\top}\over\rightarrow
\mbox{Tor}^{A\otimes A}(M\otimes M,N\otimes N)
\buildrel{{\rm Tor}^{\mu_A}(\mu_M,\mu_N)}\over\longrightarrow
\mbox{Tor}^{A}(M,N)
$$
where $\top$ is the $\top$ product (\cite[VIII.2.1]{MacLane:hom} or
\cite[XI.Proposition 1.2.1]{Cartan-Eilenberg:homalg}),
defines a natural structure of commutative graded algebra on
$\mbox{Tor}^{A}(M,N)$(\cite[Theorem 2.2]{MacLane:hom} or
\cite[XI.4 $\pitchfork$ product]{Cartan-Eilenberg:homalg}).
\begin{propriete}\label{unicite de la fibre}{\rm
\cite[A.3]{Halperin:unieal}}\qua
Suppose given a commutative diagram of augmented commutative
differential graded algebras
\[
\begin{diagram}
\node{A} \arrow{e,tb}{\varphi}{\simeq}
\arrow{s,l}{f} \node{A'} \arrow{s,r}{g} \\
\node{M} \arrow{e,tb}{\Psi}{\simeq} \node{M'}
\end{diagram}
\]
where $\varphi$ and $\Psi$ are quasi-isomorphisms.
Then
$\mbox{Tor}^{\varphi}(\Psi,{\Bbbk}):\mbox{Tor}^{A}(M,{\Bbbk})\longrightarrow\mbox{Tor}^{A'}(M',{\Bbbk});$
is an isomorphism of graded commutative algebras.
\end{propriete}
\begin{propriete}{\rm\cite[VIII.2.3]{MacLane:hom}}\qua
\label{tor d'algebres commutatives}
Consider a factorization $f=p\circ i$ where $i:A\rightarrowtail P$
is a morphism of commutative differential graded algebras such
that $P$ is an $A$-semifree
module
and $p:P\buildrel{\simeq}\over\rightarrow M$ is a quasi-isomorphism
of commutative differential graded algebras.
The homology of the commutative differential graded algebra $P\otimes_A N$,
$H_*(P\otimes_A N)$, is the graded commutative algebra $\mbox{Tor}^{A}(M,N)$.
\end{propriete}

Using this Property, Theorem A given in the Introduction
derives from the following proposition:

Let $r\geq 1$ be a integer.
Let $p$ be the characteristic of the field
${\Bbbk}$ (except when the characteristic is $0$: in this case,
we set $p=+\infty$).
We suppose now $p\not= 2$.
\begin{defin}\cite{Halperin:unieal}\qua
A topological space $X$ is {\it $(r,p)$-mild} or in the {\it Anick range}
if it is $r$-connected and 
 its homology with coefficient in ${\Bbbk}$
is concentrated in degrees $\leq rp$
and of finite type.
\end{defin}
\begin{proposition}\label{la fibre du model}
Let $f:E\rightarrow B$ be a continuous map between two topological spaces
both $(r,p)$-mild with $H_{rp}(f)$ injective.
Consider the homotopy fiber $F$ and the induced
fibration $p_0:F\twoheadrightarrow E$.
Then there are two morphisms of augmented commutative cochain algebras, denoted $A(f):A(B)\rightarrow A(E)$
and $A(p_0):A(E)\rightarrow A(F)$
such that:
\begin{enumerate}
\item There is a commutative diagram of cochain complexes
\[
\begin{diagram}
\node{C^{*}(B)}\arrow{s,r}{\simeq}\arrow{e,t}{C^{*}(f)}
\node{C^{*}(E)}\arrow{s,r}{\simeq}\arrow{e,t}{C^{*}(p_0)}
\node{C^{*}(F)}\arrow{s,r}{\simeq}\\
\node{D_1(B)}\arrow{e}
\node{D_1(E)}\arrow{e}
\node{D_1(F)}\\
\node{D_2(B)}\arrow{s,r}{\simeq}\arrow{n,r}{\simeq}\arrow{e}
\node{D_2(E)}\arrow{s,r}{\simeq}\arrow{n,r}{\simeq}\arrow{e}
\node{D_2(F)}\arrow{s,lr}{\Theta}{\simeq}\arrow{n,r}{\simeq}\\
\node{A(B)}\arrow{e,t}{A(f)}
\node{A(E)}\arrow{e,t}{A(p_0)}
\node{A(F)}
\end{diagram}
\]
where all the vertical maps are quasi-isomorphisms and where all the
maps are morphisms of augmented cochain algebras
except $\Theta:D_2(F)\buildrel{\simeq}\over\rightarrow A(F)$
who induces a morphism of graded algebras only in homology.
\item For any factorization $A(f)=\Phi\circ i$ where $i:A(B)\rightarrowtail C$ is
a morphism of augmented commutative cochain algebras such that $C$ is an $A(B)$-semifree module
and where $\Phi :C\buildrel {\simeq}\over\rightarrow A(E)$ is a quasi-isomorphism
of augmented commutative cochain algebras, there is a commutative diagram of
augmented
commutative cochain algebras
\[
\begin{diagram}
\node{A(B)}\arrow{e,t}{A(f)}\arrow{se,b}{i}
\node{A(E)}\arrow{e,t}{A(p_0)}\node{A(F)}\\
\node{}\node{C}\arrow{n,lr}{\Phi}{\simeq}\arrow{e}\arrow{se}
\node{D_3}\arrow{n,r}{\simeq}\arrow{s,r}{\simeq}\\
\node{}\node[2]{{\Bbbk}\otimes_{A(B)} C}
\end{diagram}
\]
In particular,
the cohomology $H^{*}(F)$, is isomorphic as graded algebras to
the cohomology $H^{*}({\Bbbk}\otimes_{A(B)} C)$.
\end{enumerate}
\end{proposition}
Over a field of characteristic zero,
part (1) was proved by Sullivan \cite{Sullivan:infct}
and part (2) is the Grivel-Thomas-Halperin theorem
\cite{Allday-Puppe:cohmtg,Felix-Halperin-Thomas:ratht}.

The hypotheses of Theorem A are necessary:
the space $B$ must be $(r,p)$-mild.
Indeed even for a path fibration
$\Omega X \hookrightarrow PX\twoheadrightarrow X$,
a commutative model of $X$ does not determine the cohomology algebra of the loop space.
The spaces $\Sigma \mathbb{CP} ^{p}$ and $S^{3}\vee ..\vee S^{2p+1}$
have the same commutative model but the cohomology algebras of
their loop spaces are not isomorphic.
The map $H_{rp}(f)$ must also be injective.
Take the same example as in section \ref{fibre d'une application suspendue et contreexemple}:
the suspension of the Hopf map
$\Sigma f:\Sigma S^{2p-1}\rightarrow\Sigma \mathbb{CP} ^{p-1}$.

Over a field of characteristic $p$, 
we can't improve Proposition \ref{la fibre du model}, by
``${\Bbbk}\otimes_{A(B)} C$ is weakly equivalent as a cochain
algebra to $C^{*}(F)$".
For example, let $X=K(\mathbb{Z},4)_{2p+3}$ be the $2p+3$ skeleton of a $K(\mathbb{Z},4)$.
The space $X$ is $(3,p)$-mild and
$C^{*}(\Omega X)$ is not weakly equivalent as
a cochain algebra to any commutative
cochain algebra.
Indeed, there exist two CW-complexes
denoted $Y$ and $K(\mathbb{Z},3)$ with the same $2p+2$ skeleton, respectively
homotopic to $\Omega X$ and $\Omega K(\mathbb{Z},4)$.
The two morphisms of topological monoids
$$\Omega (Y_{2p+2})\rightarrow \Omega Y\quad\mbox{and}\quad
\Omega (K(\mathbb{Z},3)_{2p+2})\rightarrow \Omega K(\mathbb{Z},3)$$ induce in homology two
algebra morphisms which are isomorphisms in degree $\leq 2p$.
Since $H_*(\Omega K(\mathbb{Z},3))\cong \Gamma \alpha _2$ as algebras,
$\Omega Y$ is $1$-connected, $H_2(\Omega Y)=\mathbb{F}_p \alpha _2$ and
$\alpha _2 ^{p}=0$.
Suppose $C^{*}(Y)$ is weakly equivalent as a cochain algebra
to a commutative cochain algebra $A$.
We can suppose that $A$ is of finite type. 
The dual of $A$, denoted $A^{\vee}$, is a cocommutative chain
coalgebra.
There is a quasi-isomorphism of chain algebras
from the cobar construction of $A^{\vee}$, denoted $\Omega(A^{\vee})$,
to $C_*(\Omega Y)$.
The Quillen construction on the coalgebra $A^{\vee}$
is a differential graded Lie algebra, denoted 
$\mathcal{L}_A$, such that $U\mathcal{L}_A:=\Omega(A^{\vee})$ \cite[p. 307 and 315]{Felix-Halperin-Thomas:ratht}.
The homology of an universal enveloping algebra of a differential
graded Lie algebra is
isomorphic as graded Hopf algebras to the universal enveloping algebra
of a graded Lie algebra~\cite[8.3]{Halperin:unieal}:
there a graded Lie algebra $E$ equipped with the following isomorphism
of graded algebras
$$H_*(\Omega Y)\cong H_*(U\mathcal{L}_A)\cong UE.$$
By the Poincar\'e-Birkoff-Witt Theorem~\cite[1.2]{Halperin:unieal},
$H_*(\Omega Y)$ admits a basis
containing $\alpha _2 ^{p}$. Thus $\alpha _2 ^{p}$ is non zero.\\
\begin{proof}[Proof of Proposition \ref{la fibre du model}]
By the naturality of Corollary~\ref{fibre} with respect to continuous maps,
we have a commutative diagram of coalgebras:
\begin{equation}\label{naturalite du theoreme de la fibre}
\begin{diagram}
\node{C_*(F)}\arrow[2]{e,t}{C_*(p_0)}
\node[2]{C_*(E)}\arrow[2]{e,t}{C_*(f)}
\node[2]{C_*(B)}\\
\node{G(f)}\arrow{n,r}{\simeq}\arrow{s,r}{\simeq}\arrow[2]{e}
\node[2]{G(E\rightarrow *)}\arrow{n,r}{\simeq}\arrow{s,r}{\simeq}\arrow[2]{e}
\node[2]{G(B\rightarrow *)}\arrow{n,r}{\simeq}\arrow{s,r}{\simeq}\\
\node{B(C_*(\Omega B);C_*(\Omega E))}\arrow[2]{e}
\node[2]{BC_*(\Omega E)}\arrow[2]{e,t}{BC_*(\Omega f)}
\node[2]{BC_*(\Omega B)}
\end{diagram}
\end{equation}
There is also a commutative diagram of augmented chain algebras \cite[Theorem~I]{Felix-Halperin-Thomas:Adamce}
\[
\begin{diagram}
\node{TX}\arrow{e,t}{\simeq}\arrow{s,l}{m(f)}
\node{\Omega C_*(E)}\arrow{s,r}{\Omega C_*(f)}\arrow{e,t}{\simeq}
\node{C_*(\Omega E)}\arrow{s,r}{C_*(\Omega f)}\\
\node{TY}\arrow{e,t}{\simeq}\node{\Omega C_*(B)}\arrow{e,t}{\simeq}
\node{C_*(\Omega B)}
\end{diagram}
\]
where $\Omega$ denotes the cobar construction,
$TX$ is a minimal
(in the sense of \cite[2.1]{Baues-Lemaire:Minmht})
free chain algebra and $m(f):TX \rightarrowtail TY$
is a minimal free extension.
Since the indecomposables functor $Q$ preserves quasi-isomorphism
between free chain algebras \cite[1.5]{Baues-Lemaire:Minmht},
$$X\cong s^{-1}H_+(E)\quad\mbox{and}\quad
Y\cong s^{-1}H_+(E)\oplus s^{-1}\mbox{coker}H_+(f)\oplus\mbox{ker}H_+(f).$$
So $X$ and $Y$ are graded vector spaces of finite type concentrated in degree $\geq r$ and $\leq rp-1$.
Denote by $m_X$ the composite $TX\buildrel{\simeq}\over\rightarrow\Omega C_*(E)
\buildrel{\simeq}\over\rightarrow C_*(\Omega E)$
and by $m_Y$ the composite $TY\buildrel{\simeq}\over\rightarrow\Omega C_*(B)
\buildrel{\simeq}\over\rightarrow C_*(\Omega B)$.
By Theorem \ref{passage aux models avec cofibration}, $m(f):TX
\rightarrowtail TY$ is an inclusion of Hopf algebras up to homotopy and
$B(m_Y;m_X):B(TY;TX)\buildrel{\simeq}\over\rightarrow B(C_*(\Omega E);C_*(\Omega B))$
is a morphism of coalgebras up to homotopy.

By Anick's Theorem~\cite[5.6]{AnickD:Hopah}, there exists a differential graded
Lie algebra $L(E)$ and an isomorphism $\varphi$ of Hopf algebras up to homotopy
between the universal envelopping algebra of $L(E)$, $UL(E)$ and $TX$.
By the naturality of Anick's Theorem with respect to Hopf algebras
up to homotopy equipped with their homotopies
(\cite{Majewski:rathmu} D.33 and D.25, see also the proof of Theorem 8.5(g)\cite{AnickD:Hopah}),
there exists a differential graded Lie algebra morphism $L(f):L(E)\rightarrow L(B)$ and a commutative
diagram of chain algebras
\[
\begin{diagram}
\node{UL(E)}\arrow{e,tb}{\cong}{\varphi}\arrow{s,l}{UL(f)}
\node{TX}\arrow{s,r}{m(f)}\\
\node{UL(B)}\arrow{e,tb}{\cong}{\Psi}\node{TY}
\end{diagram}
\]
where $\varphi$ and $\Psi$ are two algebra isomorphisms equipped
with two homotopies of algebras
$$h_{top}:(\varphi\otimes\varphi)\Delta_{UL(E)}\thickapprox_a\Delta_{TX}\varphi\quad\mbox{and}\quad
h_{bottom}:(\Psi\otimes\Psi)\Delta_{UL(B)}\thickapprox_a\Delta_{TY}\Psi$$
$$\mbox{such that}\quad h_{bottom}UL(f)=(m(f)\otimes m(f))h_{top}$$
({\it the horizontal arrows commute with the diagonals
up to natural homotopies}).
By Lemma \ref{Hopf algebras a homotopie pres}(ii), the isomorphism
$$B(\Psi;\varphi):B(UL(B);UL(E))\buildrel {\cong}\over\rightarrow B(TY;TX)$$
commutes up to chain homotopy with the diagonals.
We give the Cartan-Chevalley-Eilenberg complex
with coefficients \cite[p. 242]{Halperin:unieal} $C_*(UL(B);L(E))$ the tensor product coalgebra structure of
$UL(B)\otimes \Gamma sL(E)$.
The injection $C_*(UL(B);L(E))\buildrel {\simeq}\over\rightarrow B(UL(B);UL(E))$
is a quasi-iso\-mor\-phism of coalgebras \cite[6.11]{Felix-Halperin-Thomas:dgait}.
By functoriality of the bar construction and of the Cartan-Chevalley-Eilenberg complex
with coefficients, finally we get the commutative diagram of coalgebras up to homotopy
\begin{equation}\label{diagramme bar to complexe de Cartan}
\begin{diagram}
\node{B(C_*(\Omega B);C_*(\Omega E))}\arrow[2]{e}
\node[2]{BC_*(\Omega E)}\arrow[2]{e,t}{BC_*(\Omega f)}
\node[2]{BC_*(\Omega B)}\\
\node{B(TY;TX)}\arrow[2]{e}\arrow{n,lr}{B(m_Y;m_X)}{\simeq}
\node[2]{B(TX)}\arrow[2]{e,t}{B(m(f))}\arrow{n,lr}{B(m_X)}{\simeq}
\node[2]{B(TY)}\arrow{n,lr}{B(m_Y)}{\simeq}\\
\node{B(UL(B);UL(E))}\arrow[2]{e}\arrow{n,lr}{B(\Psi;\varphi)}{\cong}
\node[2]{B(UL(E))}\arrow[2]{e,t}{B(UL(f))}\arrow{n,lr}{B(\varphi)}{\cong}
\node[2]{B(UL(B))}\arrow{n,lr}{B(\Psi)}{\cong}\\
\node{C_*(UL(B);L(E))}\arrow[2]{e}\arrow{n,r}{\simeq}
\node[2]{C_*L(E)}\arrow[2]{e,t}{C_*L(f)}\arrow{n,r}{\simeq}
\node[2]{C_*L(B)}\arrow{n,r}{\simeq}
\end{diagram}
\end{equation}
where all the coalgebras up to homotopy are counitary and coassociative exactly except
$B(TY;TX)$, where all the morphisms commute exactly with the diagonals except $B(m_Y;m_X)$
and $B(\Psi;\varphi)$, and
where all the vertical maps are quasi-isomorphisms.
Define $A(f)$ to be $C^{*}L(f):C^{*}L(B)\rightarrow C^{*}L(E)$
and $A(p_0)$ to be the inclusion $C^{*}L(E)\hookrightarrow C^{*}(UL(B);L(E)$.
By dualizing diagram \ref{naturalite du theoreme de la fibre} and diagram
\ref{diagramme bar to complexe de Cartan}, we obtain the diagram of 1.

By the universal property of push out, there is a morphism of
commutative cochain algebras
$$C^{*}(UL(B);L(B))\otimes_{C^{*}(L(B))}C^{*}(L(E))
\buildrel{\cong}\over\longrightarrow C^{*}(UL(B);L(E))$$
which is an isomorphism since $L(B)$ is of finite type.
Recall that $C\buildrel{\simeq}\over\rightarrow C^{*}L(E)$
is a quasi-isomorphism of $C^{*}L(B)$-cochain algebras and
that $C^{*}(UL(B);L(B))$ is $C^{*}(L(B))$-semifree.
Set $D_3=C^{*}(UL(B);L(B))\otimes_{C^{*}L(B)} C$.
Then we obtain the quasi-isomorphism
$D_3\buildrel{\simeq}\over\rightarrow
C^{*}(UL(B);L(B))\otimes_{C^{*}(L(B))}C^{*}(L(E))$.
Symetrically, recall that $C^{*}(UL(B);L(B))
\buildrel{\simeq}\over\rightarrow {\Bbbk}$ is 
a quasi-isomorphism of $C^{*}L(B)$-cochain
algebras~\cite[6.10]{Felix-Halperin-Thomas:dgait} and that
 $C$ is $C^{*}(L(B))$-semifree.
Then we obtain the quasi-isomorphism $D_3\buildrel{\simeq}\over\rightarrow
{\Bbbk}\otimes_{C^{*}L(B)} C$.\end{proof}
\section{Sullivan models mod $p$}
We want to use our Theorem A for practical computations.
Like in Rational Homotopy, we need two steps.
First, we replace $A(f):A(B)\rightarrow A(E)$ by a morphism between
Sullivan models.
Second, we construct a factorization of this morphism between Sullivan
models.

Contrary to the rational case~\cite[Proposition
14.6]{Felix-Halperin-Thomas:ratht}, modulo $p$, there is in general,
no lifting lemma. Nevertheless, we have the following:
\begin{cor}\label{passage au modeles de Sullivan}
$\bullet$ Let $A(f):A(B)\rightarrow A(E)$ be a morphism of
commutative cochain algebras as in Proposition \ref{la fibre du model}.
Let $\Lambda Y$ be a Sullivan model of $A(B)$, $\Lambda X$ a Sullivan model of $A(E)$.
Then there is an acyclic commutative cochain algebra $U$ and a commutative
diagram of commutative cochain algebras
\[
\begin{diagram}
\node{}\node{\Lambda Y}\arrow{sw,t}{\Psi}\arrow{s}\arrow{e,t}{\simeq}
\node{A(B)}\arrow{s,r}{A(f)}\\
\node{\Lambda X}\node{\Lambda X\otimes U}\arrow{w,t}{\simeq}\arrow{e,t}{\simeq}
\node{A(E)}
\end{diagram}
\]
$\bullet$ Let $\Lambda Y\rightarrowtail C\buildrel {\simeq}\over\rightarrow \Lambda X$
be a factorization of $\Psi:\Lambda Y\rightarrow \Lambda X$ such that $C$
is a $\Lambda Y$-semifree module.
Then the algebra $H^{*}(F)$ is isomorphic to $H^{*}({\Bbbk}\otimes_{\Lambda Y} C)$.
(This isomorphism identifies in homology
$C^{*}(p_0):C^{*}(E)\rightarrow C^{*}(F)$ and
the quotient map $C\twoheadrightarrow {\Bbbk}\otimes_{\Lambda Y} C$.)
\end{cor}
\proof
Since $A(E)$ is concentrated in degrees $\geq r+1$ and $H^{\geq
  (r+1)p}(E)=0$,
\cite[Proposition 7.7 and Remark 7.8]{Halperin:unieal} gives
the first part of this Corollary.
For the second part, using Proposition~\ref{la fibre du model},
Property~\ref{unicite de la fibre} twice and finally Property~\ref{tor
  d'algebres commutatives}, we obtain the 
sequence of isomorphisms of graded algebras:
\begin{align}
H^{*}(F)\cong \mbox{Tor}^{A(B)}(A(E),{\Bbbk})
\cong \mbox{Tor}^{\Lambda Y}(\Lambda X\otimes U,{\Bbbk})&\notag\\
\cong \mbox{Tor}^{\Lambda Y}(\Lambda X,{\Bbbk})
&\cong H^{*}({\Bbbk}\otimes_{\Lambda Y} C).\tag*{\qed}
\end{align}

As in the rational case, we can take a factorization of $\Psi$ with relative
Sullivan models. But mod $p$, since the $p^{th}$ power of an element of even
degree is always a cycle, our relative Sullivan model will have infinitely
many generators.
We'd rather use a free divided powers algebra $\Gamma V$ where for $v\in V_{even}$,
$v^{p}=0$.
But now arises the problem of constructing morphisms of
commutative algebras from a free divided
power algebra to any commutative algebra
where the $p^{th}$ powers are not zero.
We give now an effective construction of a factorization of $\Psi$ with
divided powers algebras.
Over $\mathbb{Q}$, this factorization will be just a factorization of $\Psi$
through a minimal relative Sullivan model.

Let $A$ be a commutative graded algebra, $V$ and $W$ two graded vector spaces.
A {\it $\Gamma$-derivation} in $A\otimes \Gamma W$ is a derivation $D$
such that $D\gamma ^{k}(w)=D(w)\gamma ^{k-1}(w)$, $k\geq 1$, $w\in W^{even}$.
Any linear map $V\oplus W\rightarrow \Lambda V\otimes \Gamma W$ of degree $k$
extends to a unique $\Gamma$-derivation over $\Lambda V\otimes \Gamma W$.
\begin{lem}\label{factorisation avec algebre a puissances divisees}
Let $\Psi:(\Lambda Y,d)\rightarrow (\Lambda X,d)$ be a morphism of
commutative cochain algebras between two minimal Sullivan models
such that $X$ and $Y$ are concentrated in degree
$\geq 2$.
Then there is an explicit factorization of $\Psi$:
$$(\Lambda Y,d)\buildrel{i}\over\rightarrowtail(\Lambda Y\otimes \Lambda\mbox{coker} \varphi
\otimes\Gamma s\mbox{ker}\varphi,D)\build\twoheadrightarrow_p^{\simeq}
(\Lambda X,d)$$
where\begin{itemize}
\item $\varphi$ is the composite $Y\hookrightarrow \Lambda Y\buildrel{\Psi}\over
\rightarrow\Lambda X \twoheadrightarrow X$ and $D$ is a $\Gamma$-derivation,
\item $i$ is an 
inclusion of augmented commutative cochain algebras such that
$(\Lambda Y\otimes \Lambda\mbox{coker} \varphi\otimes\Gamma s\mbox{ker}\varphi,D)$
is $(\Lambda Y,d)$-semifree, and
\item $p$ is a surjective quasi-isomorphism of commutative cochain algebras
vanishing on $\Gamma s\mbox{ker}\varphi$.
\end{itemize}
\end{lem}
\begin{proof}We will define $$p:=\lim_{\longrightarrow}p_n.$$
We proceed by induction on $n\in \mathbb{N}^{*}$ to construct each $p_n$.
Suppose we have constructed the factorization:
$$\left(\Lambda (Y^{\leq n}),d\right)\rightarrowtail
\left(\Lambda (Y^{\leq n})\otimes\Lambda(\mbox{coker} \varphi
^{\leq n})\otimes\Gamma s(\mbox{ker}\varphi ^{\leq n}),D\right)
\build\twoheadrightarrow_{p_n}^{\simeq}
\left(\Lambda(X^{\leq n}),d\right)$$
We define now $p_{n+1}$ extending $\Psi$ and $p_n$.

Let $w\in \mbox{coker}\varphi ^{n+1}$.
Define $p_{n+1}$ in $\mbox{coker}\varphi ^{n+1}$ so that
$p_{n+1}(w)\in X^{n+1}$ represents $w$.
Then $dp_{n+1}(w)$ is a cycle of $\Lambda X^{\leq n}$.
Since $p_n$ is a surjective quasi-isomorphism, there is a cycle
$z\in \Lambda (Y^{\leq n})\otimes \Lambda(\mbox{coker} \varphi
^{\leq n})\otimes\Gamma s(\mbox{ker}\varphi ^{\leq n})$ such that $p_n(z)=dp_{n+1}(w)$. Define
$Dw=z$.

Let $v\in \mbox{ker}\varphi ^{n+1}$. Since $p_{n+1}:\Lambda (Y^{\leq
  n+1}\oplus\mbox{coker}\varphi ^{\leq n+1})\twoheadrightarrow\Lambda
  (X^{\leq n+1})$ is a surjective
morphism of graded algebras, there is
$u\in\Lambda ^{\geq 2}(Y^{\leq n}\oplus \mbox{coker}\varphi ^{\leq n})$
such that $p_{n+1}(v+u)=0$. Since $D(v+u)$ is a cycle of
$\Lambda (Y^{\leq n})\otimes \Lambda(\mbox{coker} \varphi
^{\leq n})\otimes\Gamma s(\mbox{ker}\varphi ^{\leq n})$ and $p_n$ is a surjective
quasi-isomorphism, there is $\alpha\in\Lambda (Y^{\leq n})\otimes \Lambda(\mbox{coker} \varphi
^{\leq n})\otimes\Gamma s(\mbox{ker}\varphi ^{\leq n})$ such that $p_n(\alpha)=0$ and
$D\alpha=D(v+u)$. Define $Dsv=v+u-\alpha$.

Now we have the commutative diagram of commutative cochain algebras:
\[
\begin{diagram}
\node{\Lambda (Y^{\leq n})\otimes\Lambda(\mbox{coker} \varphi
^{\leq n})\otimes\Gamma s(\mbox{ker}\varphi ^{\leq n}),D}
\arrow{s}\arrow[3]{e,tb}{p_n}{\simeq}
\node[3]{\Lambda(X^{\leq n}),d}\arrow{s}\\
\node{\Lambda (Y^{\leq n+1})\otimes\Lambda(\mbox{coker} \varphi
^{\leq n+1})\otimes\Gamma s(\mbox{ker}\varphi ^{\leq n+1}),D}
\arrow{s}\arrow[3]{e,t}{p_{n+1}}
\node[3]{\Lambda(X^{\leq n+1}),d}\arrow{s}\\
\node{\Lambda (Y^{n+1})\otimes \Lambda(\mbox{coker} \varphi
^{n+1})\otimes\Gamma s(\mbox{ker}\varphi ^{n+1}),\overline{D}}
\arrow[3]{e,tb}{\overline{p_{n+1}}}{\simeq}
\node[3]{\Lambda (X^{n+1}),0}
\end{diagram}
\]
Since $p_n$ and $\overline{p_{n+1}}$ are quasi-isomorphisms, by comparison
of the $E_2$-term of the algebraic Serre spectral sequence
associated to each column, $p_{n+1}$ is
a quasi-isomorphism.\end{proof}
\begin{ex}
Let $f:S^{2}\hookrightarrow \mathbb{CP} ^{n}$ be the inclusion of CW-complexes with $n\geq 2$.
Applying Corollary \ref{passage au modeles de Sullivan}, $\Psi$ is
$(\Lambda(x_2,y_{2n+1}),d)\rightarrow
(\Lambda(x_2,z_3),d)$ with $dy_{2n+1}=x_2^{n+1}$ and $dz_3=x_2^{2}$.
Thus $\Psi y_{2n+1}=z_3x_2^{n-1}$.
By Lemma \ref{factorisation avec algebre a puissances divisees},
$\Psi$ factors through the commutative cochain algebra
$(\Lambda(x_2,y_{2n+1},z_3)\otimes\Gamma sy_{2n+1},D)$ with $Dz_3=x_2^{2}$ and
$Dsy_{2n+1}=y_{2n+1}-z_3x_2^{n-1}$.
So $H^{*}(F)\cong\Lambda z_3\otimes\Gamma sy_{2n+1}$ for $p\geq 2n$.
\end{ex}
\section{Proof of Theorem B}
The key to the proof of Theorem A is to apply
Anick's Theorem~\cite[5.6]{AnickD:Hopah}.
One of the goals of Anick for developing this theorem was to prove
a result suggested by McGibbon and Wilkerson~\cite[p. 699]{McGibbon-Wilkerson:loosfclp}:
``If $X$ is a finite simply-connected CW-complex then for large primes,
$p^{th}$ powers vanish in
$\tilde{H}^{*}(\Omega X;\mathbb{F}_p)$.''
Anick~\cite[9.1]{AnickD:Hopah} proved precisely that
``If $X$ is $(r,p)$-mild then
$p^{th}$ powers vanish in  $\tilde{H}^{*}(\Omega X;\mathbb{F}_p)$.''.
After Anick, Halperin proved
in~\cite[Theorem 8.3 and Poincar\'e-Birkoff-Witt Theorem]{Halperin:unieal} 
that in fact:
\begin{cor}{\rm\cite{Halperin:unieal}}\qua
If $X$ is $(r,p)$-mild then
the algebra $H^{*}(\Omega X)$ is isomorphic to $\Gamma sV$
where $\Lambda V$ is a minimal Sullivan model of $A(X)$.
\end{cor}
\begin{proof} Apply Corollary \ref{passage au modeles de Sullivan} to
$*\rightarrow X$.
Consider the factorization of $(\Lambda V,d)\twoheadrightarrow
({\Bbbk},0)$,
$$(\Lambda V,d)\rightarrowtail
(\Lambda V\otimes\Gamma sV,D)\buildrel{\simeq}
\over\twoheadrightarrow {\Bbbk}$$
given by Lemma \ref{factorisation avec algebre a puissances divisees}.
See that the cofiber
$({\Bbbk},0)\otimes_{(\Lambda V,d)} (\Lambda V\otimes \Gamma sV,D)$ has a null
differential \cite[2.6]{Halperin:unieal}.\end{proof}

Actually, we can show now that Anick's result on $p^{th}$ powers and
Halperin's result on a divided powers algebra structure remain
valid if we consider the fiber of any fibration in the Anick range instead
of just the loop fibration.
But before we need some definitions concerning divided powers algebras
with differential.

A {\it differential divided powers algebra} or {\it $\Gamma$-algebra}
is a 
commutative cochain algebra $A$ equipped with a system
$(\gamma^{k})_{k\in\mathbb{N}}$
of divided powers \cite[page 124]{Brown:cohgro} such that
$d\gamma^{k}(a)=d(a)\,\gamma^{k-1}(a)$.
Let $A$, $B$ be two $\Gamma$-algebras.
A {\it $\Gamma$-morphism} $f:A\rightarrow B$ is
a morphism of augmented commutative cochain algebras such that
$f\gamma^{k}(a)=\gamma^{k}f(a)$.

A {\it $\Gamma$-free extension} is an inclusion of augmented
commutative cochain algebras: $(A,d)\rightarrowtail(A\otimes \Gamma V,D)$
such that $\displaystyle V=\oplus_{k\in \mathbb{N}}V(k)$,
$D:V(k)\rightarrow A\otimes \Gamma V(<k),\:k\in\mathbb{N}$
and $D$ is a $\Gamma$-derivation.
In particular, if $A$ is a $\Gamma$-algebra,
than a $\Gamma$-free extension
$(A,d)\rightarrowtail(A\otimes \Gamma V,D)$
is a $\Gamma$-morphism.

A commutative cochain algebra (respectively $\Gamma$-algebra) $A$
is {\it admissible} (respectively {\it $\Gamma$-admissible})
if there is a surjective morphism of commutative co\-chain algebras
(respectively $\Gamma$-morphism) $C\twoheadrightarrow A$
with $C$ acyclic.
\begin{propriete}\label{surjection entre admissibles}{\rm\cite[II.2.6]{Halperin:divpa}}\qua
Let $f:A\rightarrow B$ be a morphism of commutative cochain algebras
(respectively a $\Gamma$-morphism). If $f$ is surjective and $A$
is admissible (respectively $\Gamma$-admissible) then so is $B$.
\end{propriete}
\begin{proposition}\label{admissible et puissances divisees}{\rm\cite[II.2.7]{Halperin:divpa}}\qua
\begin{itemize}
\item[\rm(i)] If $f:A\rightarrow B$ is a morphism of commutative
cochain algebras with $B$ admissible
then we have the commutative diagram of commutative cochain algebras
\[
\begin{diagram}
\node{A}\arrow[2]{e}\arrow{ese}\arrow{s}\node[2]{B}\\
\node{A\otimes\Gamma V'}
\node[2]{A\otimes\Lambda V}\arrow[2]{w,b}{\simeq}\arrow{n,r}{\simeq}
\end{diagram}
\]
where $A\rightarrowtail A\otimes\Lambda V$ is a relative
Sullivan model and
$A\rightarrowtail A\otimes\Gamma V'$ is a
$\Gamma$-free extension.
\item[\rm(ii)] In particular, if $B$ is any admissible commutative cochain
  algebra, there are quasi-isomorphisms of commutative cochain algebras
$$\Gamma V'\buildrel{\simeq}\over\longleftarrow\Lambda V
\buildrel{\simeq}\over\longrightarrow B$$
where $\Gamma V'$ is a $\Gamma$-algebra.
\end{itemize}
\end{proposition}

The essential role of $\Gamma$-admissible algebras
is that
\begin{propriete}\label{homology admissible}{\rm\cite[1.3]{{Avramov-Halperin:lookingglass}}}\qua
If $A$ is a $\Gamma$-admissible algebra
then $H(A)$ is a divided powers algebra (not true
if $A$ was only a $\Gamma$-algebra!).
\end{propriete}
\begin{lem}\label{condition d'admissibilite}
Let $A$ be a commutative cochain algebra.
Assume that for some $r\geq 1$, A satisfies
$A={\Bbbk}\oplus\{A^{i}\}_{i\geq r}$.
\begin{itemize}
\item[\rm(i)] If $H^{i}(A)=0,\: i\geq rp$,
then $A$ is admissible.
\item[\rm(ii)] If $A$ is a $\Gamma$-algebra
and $H^{i}(A)=0,\: i\geq rp+p-1$,
then $A$ is $\Gamma$-admissible.
\end{itemize}
\end{lem}
\begin{proof} (i)\qua This lemma is just a slight improvement
from \cite[Lemma 7.6]{Halperin:unieal}
and the proof is the same.

(ii)\qua After replacing free commutative algebras by free divided powers algebras,
the proof is the same as in (i).\end{proof}
\begin{lem}\label{Tor algebre a puissances divisees}
Let $A$ and $M$ be two commutative cochain algebras concentrated in degrees $\geq r+1$.
Consider a morphism of algebras $A\rightarrow M$.
If $H^{\geq rp+p}(A)=H^{\geq rp+p-1}(M)=0$ then
$\mbox{Tor}^{A}(M,{\Bbbk})$ is a divided powers algebra.
\end{lem}
\begin{proof} By Lemma \ref{condition d'admissibilite} (i),
$A$ and $M$ are admissible.
By Proposition \ref{admissible et puissances divisees} (ii),
there are quasi-isomorphisms of commutative cochain algebras
$$\Gamma X'\buildrel{\simeq}\over\longleftarrow\Lambda X
\buildrel{\simeq}\over\longrightarrow A$$
where $X$ and $X'$ are concentrated in degrees $\geq r+1$.
By Proposition \ref{admissible et puissances divisees} (i),
we get the commutative diagram of commutative cochain algebras
\[
\begin{diagram}
\node{A}\arrow{e}\node{M}\\
\node{\Lambda X}\arrow{n,l}{\simeq}\arrow{e}\arrow{se}
\node{\Lambda X\otimes\Lambda Y}\arrow{n,r}{\simeq}
\arrow{s,r}{\simeq}\\
\node{}\node{\Lambda X\otimes\Gamma Y'}
\end{diagram}
\]
where $Y$ and $Y'$ are concentrated in degrees $\geq r$.
Since $\Lambda X\rightarrowtail\Lambda X\otimes\Gamma Y'$ is a
$\Gamma$-free extension,
$\Lambda X\otimes\Gamma Y'$ is $\Lambda X$-semifree.
Therefore, by push-out, we have the commutative diagram of commutative
cochain algebras
\[
\begin{diagram}
\node{\Lambda X}\arrow[2]{e}\arrow{s,l}{\simeq}
\node[2]{\Lambda X\otimes\Gamma Y'}\arrow{s,r}{\simeq}\\
\node{\Gamma X'}\arrow[2]{e}\node[2]{\Gamma X'\otimes\Gamma Y'}
\end{diagram}
\]
where $\Lambda X\otimes\Gamma Y'
\buildrel{\simeq}\over\longrightarrow\Gamma X'\otimes\Gamma Y'$ is
a quasi-isomorphism~\cite[2.3(i)]{Felix-Halperin-Thomas:dgait}.
Since push-outs preserve $\Gamma$-free extension,
$\Gamma X'\rightarrowtail\Gamma X'\otimes\Gamma Y'$
is a $\Gamma$-free extension.
So $\Gamma X'\otimes\Gamma Y'$ is $\Gamma X'$-semifree,
and by Property \ref{unicite de la fibre}, the cohomology
algebra of the cofiber $\Gamma Y'$ is $\mbox{Tor}^{A}(M,{\Bbbk})$.
Now since $\Gamma X'$ is a $\Gamma$-algebra,
so is $\Gamma X'\otimes\Gamma Y'$.
Since $\Gamma X'\otimes\Gamma Y'$ is concentrated in
degrees $\geq r$ and its cohomology is null in degrees
$\geq rp+p-1$, by Lemma \ref{condition d'admissibilite}(ii),
$\Gamma X'\otimes\Gamma Y'$ is $\Gamma$-admissible.
Since $\Gamma X'\otimes\Gamma Y'\twoheadrightarrow
{\Bbbk}\otimes_{\Gamma X'}(\Gamma X'\otimes\Gamma Y')=\Gamma Y'$
is a surjective $\Gamma$-morphism,
by Property \ref{surjection entre admissibles},
$\Gamma Y'$ is a $\Gamma$-admissible.
So by Property \ref{homology admissible},
$H(\Gamma Y')$ is a divided powers algebra.\end{proof}
\begin{theor}\label{relative Wilkerson's conjecture}
Let $p$ be an odd prime and
let $f:E\twoheadrightarrow B$ be a fibration of fiber $F$ such that
$E$ and $B$ are
both $(r,p)$-mild with $H_{rp}(f)$ injective.
Then the cohomology algebra $H^{*}(F;\mathbb{F}_p)$
is a (not necessarily free!) divided powers algebra.
In particular, $p^{th}$ powers vanish in the reduced cohomology
$\tilde{H}^{*}(F;\mathbb{F}_p)$.
\end{theor}
\begin{proof} By Theorem A,
$H^{*}(F;{\Bbbk})\cong \mbox{Tor}^{A(B)}(A(E),{\Bbbk})$.
Since $A(B)$ and $A(E)$ are concentrated in degrees $\geq r+1$
and their cohomology is null in degrees $\geq rp$,
by Lemma \ref{Tor algebre a puissances divisees},
$\mbox{Tor}^{A(B)}(A(E),{\Bbbk})$ is a divided powers algebra.\end{proof}
\providecommand{\bysame}{\leavevmode\hbox to3em{\hrulefill}\thinspace}

\Addresses
\recd

\end{document}